\documentstyle[11pt]{article}
\newfam\msbfam
\font\tenmsb=msbm10    \textfont\msbfam=\tenmsb \font\sevenmsb=msbm7
\scriptfont\msbfam=\sevenmsb \font\fivemsb=msbm5
\scriptscriptfont\msbfam=\fivemsb
\def\Bbb{\fam\msbfam \tenmsb}

\def\rr{{\Bbb R}}
\def\rz{{{\rr}^n}}
\def\zz{{\Bbb Z}}
\def\nn{{\Bbb N}}
\def\cc{{\Bbb C}}

\def\fz{\infty}
\def\az{\alpha}
\def\supp{{\rm{\ supp\ }}}

\def\ez{\epsilon}

\def\pz{\partial}
\def\tz{\theta}

\def\vz{\varphi}
\def\dz{\delta}
\def\gz{\gamma}

\def\lz{\lambda}
\def\oz{\Omega}
\def\supp{{\rm supp}}

\def\wt{\widetilde}
\def\wz{\omega}
\def\l{\left}
\def\r{\right}

\def\dsum{\displaystyle\sum}

\def\dint{\displaystyle\int}
\def\dfrac{\displaystyle\frac}
\def\dsup{\displaystyle\sup}
\def\dlim{\displaystyle\lim}
\def\dinf{\displaystyle\inf}

\newtheorem{thm}{\hskip\parindent Theorem}

\newtheorem{lem}{\hskip\parindent Lemma}

\newtheorem{prop}{\hskip\parindent Proposition}

\begin{document}

\baselineskip=15pt
\renewcommand{\arraystretch}{2}
\arraycolsep=1.2pt

\title{ Weighted norm inequalities,  spectral multipliers and Littlewood-Paley operators in the  Schr\"odinger settings
 \footnotetext{ \hspace{-0.65
cm} 2000 Mathematics Subject  Classification: 42B25, 42B20.\\
The  research was supported  by the NNSF (10971002) of China.\\}}
%The project was supported by  the NNSF(10371004) of China.}}

\author{ Lin Tang}
\date{}
\maketitle

{\bf Abstract}\quad   In this paper, we establish a good-$\lz$ inequality with two parameters in the Schr\"odinger settings. As it's applications, we obtain weighted estimates for spectral multipliers  and Littlewood-Paley operators and their commutators in the Schr\"odinger settings.
\bigskip

\begin{center}{\bf 1. Introduction }\end{center}

In this paper, we consider the  divergence Schr\"odinger differential operator
$$L=-\pz_i(a_{ij}(x)\pz_j)+V(x)\ {\rm on}\ \rz,\ \ n\ge 3,$$ where  $V(x)$
is a  nonnegative potential satisfying  certain reverse H\"older
class.   In this paper, we always assume that the coefficients of these operators are bounded and measurable, and $a_{ij}$ are real symmetric, uniformly elliptic, i.e., for some $\dz\in(0,1]$,
$$ a_{ij}=a_{ji},\ |a_{ij}|\le\dz^{-1},\ \dz|\xi|^2\le a_{ij}\xi_i\xi_j\le\dz^{-1}|\xi|^2. \eqno(1.1)$$

 We say a nonnegative locally $L^q$
integral function $V(x)$ on $\rr^n$ is said to belong to $B_q(1<q\le
\fz)$ if there exists $C>0$ such that the reverse H\"older
inequality
$$\l(\dfrac 1{|B(x,r)|}\dint_{B(x,r)} V^q(y)dy\r)^{1/q}\le C\l(\dfrac
1{|B(x,r)|}\dint_{B(x,r)} V(y)dy\r)$$ holds for every $x\in\rz$ and
$0<r<\fz$, where $B(x,r)$ denotes the ball centered at $x$ with
radius $r$. In particular, if $V$ is a nonnegative polynomial, then
$V\in B_\fz$.  Throughout this paper, we always assume that
$0\not\equiv V\in B_n/2$.

The study of Schr\"odinger operator $L_0=-\triangle+V$ recently
attracted much attention; see \cite{aa,bhs1,bhs2,dz1,s1,yz,z}.
In particular, it should be pointed out that Shen \cite{s1} proved
the Schr\"odinger type operators, such as
$\nabla(-\Delta+V)^{-1}\nabla$, $\nabla(-\Delta+V)^{-1/2}$,
$(-\Delta+V)^{-1/2}\nabla$  with $V\in B_n$, $(-\Delta+V)^{i\gz}$  with $\gz\in\rr$ and $V\in B_{n/2}$,  are standard
Calder\'on-Zygmund operators. Later, Auscher and Ali \cite{aa} improved some results of Shen \cite{s1}.

 Recently, Bongioanni, etc,  \cite{bhs1} proved
$L^p(\rz)(1<p<\fz)$ boundedness for commutators of Riesz transforms
associated with Schr\"odinger operator with
$BMO_\tz(\rho)$ functions which include the class $BMO$ function,
and they \cite{ bhs2} established the weighted boundedness
for
Riesz transforms, fractional integrals  and Littlewood-Paley functions
 associated with Schr\"odinger operators with
weight $A_p^{\rho,\tz}$ class which includes the Muckenhoupt weight class.
Very recently, Tang \cite{t1,t2} established a new good-$\lz$ inequality, and
 obtained the weighted norm inequalities for some Schr\"odinger type operators, which include
commutators for  Riesz transforms, fractional integrals  and Littlewood-Paley functions
 associated with Schr\"odinger operators.

It should be pointed out that the results above were obtained by using sizes estimates of kernels of Schr\"odinger type operators. In some cases, we may meet some
Schr\"odinger operators that they do not have an integral representation by a kernel with sizes estimates. To deal with latter case, in this paper,
we will establish a good-$\lz$ inequality with two parameters  in the Schr\"odinger settings. As it's applications, we obtain weighted estimates for spectral multipliers and Littlewood-Paley operators and their commutators in the Schr\"odinger settings.

The paper is organized as follows. In Section 2,  we give some notation and basic results.  In Section 3, we establish a good-$\lz$ inequality with two parameters in the Schr\"odinger settings. In Section 4,  we obtain
weighted  inequalities for spectral multipliers and their commutators in the Schr\"odinger settings. Finally, we give
weighted  estimates for Littlewood-Paley operators and their commutators in the Schr\"odinger settings. in Section 5.

Throughout this paper, we let $C$ denote  constants that are
independent of the main parameters involved but whose value may
differ from line to line. By $A\sim B$, we mean that there exists a
constant $C>1$ such that $1/C\le A/B\le C$.

\bigskip

\begin{center} {\bf 2. Preliminaries}\end{center}

We first recall some notation.  Given $B=B(x,r)$
and $\lz>0$, we will write $\lz B$ for the $\lz$-dilate ball,
which is the ball with the same center $x$ and with radius $\lz
r$. Similarly, $Q(x,r)$ denotes the cube centered at $x$ with the
sidelength $r$ (here and below only cubes
with sides parallel to the coordinate axes are considered), and $\lz Q(x,r)=Q(x,\lz r)$.
 Given a Lebesgue
measurable set $E$ and a weight $\wz$, $|E|$ will denote the
Lebesgue measure of $E$ and $\wz(E)=\int_E\wz dx$. For $0< p<\fz$
$$\|f\|_{L^p(\wz)}=\l(\int_\rz |f(y)|^p\wz(y)dy\r)^{1/p},\ \|f\|_{L^{p,\fz}(\wz)}=\dsup_{\lz>0}\lz^{1/p}\wz\{ |f(y)|>\lz\}.$$
If $\wz=1$, we simply denote $\|f\|_{L^p}=\|f\|_{L^p(\wz)},\ \|f\|_{L^{p,\fz}}=\|f\|_{L^{p,\fz}(\wz)}$.

 The function $m_V(x)$ is defined by
$$\rho(x)=\dfrac 1{m_V(x)}=\dsup_{r>0}\l\{r:\ \dfrac
{1}{r^{n-2}}\dint_{B(x,r)}V(y)dy\le 1\r\}.$$Obviously,
$0<m_V(x)<\fz$ if $V\not=0$. In particular, $m_V(x)=1$ with $V=1$
and $m_V(x)\sim (1+|x|)$ with $V=|x|^2$.

\begin{lem}\label{l2.1.}\hspace{-0.1cm}{\rm\bf 2.1(\cite{s1}).}\quad
There exists $l_0>0$ and $C_0>1$such that $$\dfrac
1{C_0}\l(1+|x-y|/\rho(x)\r)^{-l_0}\le \dfrac{\rho(y)}{\rho(x)}\le
C_0\l(1+|x-y|/\rho(x)\r)^{l_0/(l_0+1)}.
$$ In particular, $\rho(x)\sim \rho(y)$ if $|x-y|<C\rho(x)$.
\end{lem}
In this paper, we write $\Psi_\tz(B)=(1+r/\rho(x_0))^\tz$, where
 $\tz> 0$, $x_0$ and $r$ denotes the center and radius of $B$ respectively.

 A weight will always mean a nonnegative function which is
locally integrable. As in \cite{bhs2}, we say that a weight $\wz$ belongs to the
class $A^{\rho,\tz}_p$ for $1<p<\fz$, if there is a constant $C$ such that
for all ball
  $B=B(x,r)$
$$\l(\dfrac 1{\Psi_\tz(B)|B|}\dint_B\wz(y)\,dy\r)
\l(\dfrac 1{\Psi_\tz(B)|B|}\dint_B\wz^{-\frac 1{p
-1}}(y)\,dy\r)^{p-1}\le C.$$
 We also
say that a  nonnegative function $\wz$ satisfies the $A^{\rho,\tz}_1$
condition if there exists a constant $C$ for all balls $B$
$$M_{\rho,\tz}(\wz)(x)\le C \wz(x), \ a.e.\ x\in\rz.$$
 where
$$M_{\rho,\tz} f(x)=\dsup_{x\in B}\dfrac 1{\Psi_\tz(B)|B|}\dint_B|f(y)|\,dy.$$
   Since $\Psi_\tz(B)\ge 1$, obviously,
$A_p\subset A_p^{\rho,\tz}$ for $1\le p<\fz$, where $A_p$ denote
the classical Muckenhoupt weights; see \cite{gr} and \cite{m}. We
will see that $A_p\subset\subset A_p^{\rho,\tz}$ for $1\le p<\fz$ in
some cases.  In fact, let $\tz>0$ and $0\le\gz\le\tz$, it is easy
to check that $\wz(x)=(1+|x|)^{-(n+\gz)}\not\in A_\fz=\bigcup_{p\ge 1}A_p$ and
$\wz(x)dx$ is not a doubling measure, but
$\wz(x)=(1+|x|)^{-(n+\gz)}\in A_1^{\rho,\tz}$ provided that  $V=1$ and
$\Psi_\tz(B(x_0,r))=(1+r)^\tz$.

When $V=0$ and $\tz=0$, we
denote $M_{0,0}f(x)$ by $Mf(x)$( the standard Hardy-Littlewood
maximal function). It is easy to see that $|f(x)|\le M_{\rho,\tz}
f(x)\le Mf(x)$ for $a.e.\ x\in\rz$ and $\tz\ge 0$.

Similar to the classical Muckenhoupt weights(see \cite{gr,s}), we give some
properties for weight class $A_p^{\rho,\tz}$ for $p\ge 1$.
\begin{prop}\label{p2.1.}\hspace{-0.1cm}{\rm\bf 2.1.}\quad
Let $\wz\in A^{\rho,\fz}_p=\bigcup_{\tz\ge
0}A_p^{\rho,\tz}$ for $p\ge 1$. Then
\begin{enumerate}
\item[(i)]If $ 1\le p_1<p_2<\fz$, then $A_{p_1}^{\rho,\tz}\subset
A_{p_2}^{\rho,\tz}$. \item[(ii)] $\wz\in A_p^{\rho,\tz}$ if and only
if $\wz^{-\frac 1{p-1}}\in A_{p'}^{\rho,\tz}$, where $1/p+1/p'=1.$
\item[(iii)] If $\wz\in A_p^{\rho,\fz},\ 1<p<\fz$, then there exists
$\ez>0$ such that $\wz\in A_{p-\ez}^{\rho,\fz}.$
\end{enumerate}
\end{prop}
{\it Proof.}\quad (i) and (ii) are obvious by the definition of
$A_p^{\rho,\tz}$. (iii) is proved in \cite{bhs2}. In fact, from
Lemma 5 in \cite{bhs2}, we know that if $\wz\in A_p^{\rho,\tz}$,
then $\wz\in A_{p_0}^{\rho,\tz_0}$, where $p_0=1+\frac
{p-1}{\dz}<p$ with some $\dz>1$($\dz$ is a constant depending only on
the $RH_{\dz}^{loc}$ constant of $\wz$, see below) and $\tz_0=\frac {\tz p+\eta
(p-1)}{p_0}$ with $\eta=\tz
p+(\tz+n)\frac{pl_0}{l_0+1}+(l_0+1)n$.\hfill$\Box$

The local reverse H\"older classes are defined in the following way: $\wz\in RH_q^{loc}$, $1<q<\fz$, if there is a constant $C$ such that for every ball $B(x_0,r) \subset\rz$ with $r<\rho(x_0)$,
$$\l(\dfrac 1{|B|}\dint_B w(x)^q dx\r)^{1/q}\le C\dfrac 1{|B|}\dint_B w(x)dx.$$
The endpoint $q=\fz$ is given by the condition: $\wz\in RH_\fz^{loc}$ whenever, for every ball $B(x_0,r) \subset\rz$ with $r<\rho(x_0)$,
$$\wz(x)\le C\dfrac 1{|B|}\dint_B w(x)dx,\quad {\rm for\ a.e.}\ x\in B.$$ From Lemma 5 in \cite{bhs2}, we know that if $\wz\in A^{\rho,\fz}_p$ for $p\ge 1$, then there exists a $q>1$, such that $\wz\in RH^{loc}_q$. In addition, it is easy to see that if  $\wz\in RH_q^{loc}$ with $1<q<\fz$, then there exists $\ez>0$ such that $\wz\in RH_{q+\ez}^{loc}$.

Next we give some weighted estimates for $M_{\rho,\tz}$.
\begin{lem}\label{l2.2.}\hspace{-0.1cm}{\rm\bf 2.2(\cite{t1}).}\quad
Let $1\le p_1 <\fz$ and suppose that $\wz\in A_{p_1}^{\rho,\tz}$. If $p_1<p<\fz$,
then the equality
$$\dint_\rz|M_{\rho,\tz} f(x)|^{p}\wz(x)dx\le
C\dint_\rz|f(x)|^{p}\wz(x)dx.$$ Furthermore, let $1\le p<\fz$,
$\wz\in A_p^\tz$   if and only if
$$\wz(\{x\in\rz:\ M_{\rho,\tz}f(x)>\lz\})\le \dfrac
{C}{\lz^p}\dint_\rz|f(x)|^p\wz(x)dx.$$
\end{lem}

\begin{lem}\label{l2.3.}\hspace{-0.1cm}{\rm\bf 2.3(\cite{t1}).}\quad
Let $1<p<\fz$, $p'=p/(p-1)$ and assume that $\wz\in A_p^{\rho,\tz}$.
 There
exists a constant $C>0$ such that
$$\|M_{\rho,p'\tz}f\|_{L^p(\wz)}\le C\|f\|_{L^p(\wz)}.$$
\end{lem}

In addition,
Bongioanni, etc, \cite{bhs1} introduce a new space $BMO_\tz(\rho)$ defined by
$$\|f\|_{BMO_\tz(\rho)}=\dsup_{B\subset
\rz}\dfrac 1{\Psi_\tz(B)|B|}\dint_B|f(x)-f_B|dx<\fz,$$ where
$f_B=\frac 1{|B|}\int_Bf(y)dy$.

In particularly, Bongioanni, etc, \cite{bhs1} proved the following result for $BMO_\tz(\rho)$.
\begin{lem}\label{l2.4.}\hspace{-0.1cm}{\rm\bf 2.4.}\quad
Let $\tz>0$ and $1\le s<\fz$. If $b\in BMO_\tz(\rho)$, then
$$\l(\dfrac 1{|B|}\dint_B|b-b_B|^s\r)^{1/s}\le C_{\tz,s}\|b\|_{BMO_\tz(\rho)}\l(1+\dfrac r{\rho(x)}\r)^{\tz'},$$ for all
$B=B(x,r)$, with
$x\in\rz$ and $r>0$, where $\tz'=(l_0+1)\tz$.
\end{lem}

Now we define $BMO_\fz(\rho)=\bigcup_{\tz>0}BMO_\tz(\rho)$. Obviously, the classical $BMO$ is properly contained in $BMO_\tz(\rho)$; more examples see \cite{bhs1}.

Applying Lemma 2.4, Tang \cite{t1} proved the following John-Nireberg inequality for $BMO_\tz(\rho)$.
\begin{prop}\label{p2.2.}\hspace{-0.1cm}{\rm\bf 2.2.}\quad
 Suppose that $f$ is in $BMO_\tz(\rho)$. There exist
positive constants $\gz$ and $C$ such that
$$\dsup_{ B}\dfrac 1{|B|}\dint_B\exp\l\{\frac
\gz{\|f\|_{BMO_\tz(\rho)}\Psi_{\tz'}(B)}|f(x)-f_B|\r\}\,dx\le C.$$
\end{prop}

We remark that  balls can be replaced by cubes  in definitions of
$A_p^{\rho,\tz}$, $BMO_\tz(\rho)$ and $M_{\rho,\tz}$, since
$\Psi_\tz(B)\le  \Psi_\tz(2B)\le 2^\tz \Psi_\tz(B)$.

Finally, we recall some basic definitions and facts about Orlicz
spaces, referring to \cite{rr} for a complete account.

 A function $B(t): [0,\fz)\to [0,\fz)$ is called a Young function if
it is continuous, convex, increasing and satisfies $\Phi(0)=0$ and
$B\to \fz$ as $t\to\fz$. If $B$ is a Young function, we define the
$B$-average of a function $f$ over a cube(ball) $Q$ by means of the
following Luxemberg norm:
$$\begin{array}{cl}
\|f\|_{B,Q,\wz}&=\dinf\l\{\lz>0:\ \dfrac
1{\wz(Q)}\dint_QB\l(\dfrac {|f(y)|}\lz\r)\wz(y)\,dy\le 1\r\}\\
&\sim \dinf\l\{t>0:\ t+\dfrac
t{\wz(Q)}\dint_QB\l(\dfrac {|f(y)|}t\r)\wz(y)\,dy\r\},\end{array}$$
where $\wz(y) dy$ is Borel measure.
 If $A,\ B$
and $C$ are Young functions such that
$$A^{-1}(t)B^{-1}(t)\le C^{-1}(t),$$
where $ A^{-1}$ is the complementary Young function associated to
$A$, then
$$\|fg\|_{C,Q,\wz}\le 2\|f\|_{A,Q,\wz}\|g\|_{B,Q,\wz}.$$
The examples to be considered in our study will be
$A^{-1}(t)=\log(1+t),\ B^{-1}(t)=t/\log(e+t)$ and $C^{-1}(t)=t$.
Then $A(t)\sim e^t$ and $B(t)\sim t\log(e+t)$, which gives the
generalized H\"older's inequality for any cubes(balls) $Q$
$$\dfrac
1{\wz(Q)}\dint_Q|fg|\wz(y)\,dy\le \|f\|_{A,Q,\wz}\|g\|_{B,Q,\wz}$$ holds for any  Borel measure $\wz(y) dy$ . And we define the corresponding maximal function
$$M_B f(x)=\dsup_{Q:x\in Q}\|f\|_{B,Q}:=\dsup_{Q:x\in Q}\|f\|_{B,Q,1}$$
and for $0<\eta<\fz$
$$M_{B,\rho,\eta} f(x)=\dsup_{Q:x\in Q}(\Psi_\eta(Q))^{-1}\|f\|_{B,Q}.$$

 The  examples such as
$B(t)=t(1+\log^+t)^\az(\az>0)$ with the maximal function denoted by
$M_{L(\log L)^\az,\rho,\eta}$. The complementary Young function is given by $\bar
B(t)\approx e^{t^\az}$ with the corresponding maximal function denoted
by $M_{\exp L^\az,\rho,\eta}$. In this previous case, it is well known that for $k\ge 1$, from the the proof of Lemma 4.1 in \cite{t1}, we have
 $$M_{L(\log L)^k,\rho,\eta}f\le C_{\eta,k}M^{k+1}_{\rho,\eta/2^k}f,\eqno(2.1)$$
where $M^{k+1}_{\rho,\eta/2^k}$ is the $k+1$-iteration of $M_{\rho,\eta/2^k}$.

For these
example and using Theorem 2.1, if $b\in BMO_{\eta}(\rho)$ and $b_Q$ denotes
its average on the cube(ball) $Q$, then
$\|(b-b_Q)/\Psi_{\eta'}(Q)\|_{expL, Q}\le C\|b\|_{BMO_{\eta}(\rho)},$
where $\eta'=(1+l_0)\eta>0$. This yields the following estimates: First, for each cube(ball) $Q$ and $x\in Q$
$$\begin{array}{cl}
\dfrac 1{\Psi(Q)^{\eta'(kp_0+1)}|Q|}&\dint_Q|b-b_Q|^{kp_0}|f|^{p_0}dy\\
&\le \|(b-b_Q)/\Psi_{\eta'}(Q)\|^{kp_0}_{expL, Q}
\||f|^{p_0}/\Psi_{\eta'}(Q)\|_{L(\log L)^{kp_0},Q}\\
&\le C\|b\|_{BMO_{\eta}(\rho)}^{kp_0}M_{L(\log L)^{kp_0},\rho,\eta'}(|f|^{p_0})(x)\\
&\le  C\|b\|_{BMO_{\eta}(\rho)}^{kp_0}M^{[kp_0]+2}_{\rho,\eta'/2^{[kp_0]+1}}(|f|^{p_0})(x).
\end{array}\eqno(2.2)$$
where $[s]$ is the integer part of $s$. Second, for $j\ge 1$ and each $Q$,
$$\begin{array}{cl}
\|(b-b_Q)/\Psi_{\eta'}(2^jQ)\|_{expL,\rho, 2^jQ}&\le \|(b-b_{2^jQ})/\Psi_{\eta'}(2^jQ)\|_{expL,\rho, 2^jQ}
+|b_{2^jQ}-b_Q|/\Psi_{\eta'}(2^jQ)\\
&\le Cj\|b\|_{BMO_\eta(\rho)}.
\end{array}\eqno(2.3)$$
 In addition,  for any  cube(ball) $Q=Q(x_0,r)$ with $r<\rho(x_0)$, let $b\in BMO_\tz(\rho)$ and $b_Q$ denotes
its average on $Q$, if $\wz\in RH_q^{loc}$ for some $q>1$, then by Proposition 2.2,
$$\|(b-b_Q)\|_{expL, Q,\wz}\le C\|b\|_{BMO_\tz(\rho)},\eqno(2.4)$$
and $$\dfrac 1{\wz(Q)}\dint_Q|b(y)-b_Q|\wz(y)dy\le C\|b\|_{BMO_\tz(\rho)}.\eqno(2.5)$$

\begin{center} {\bf 3. Two-parameter good-$\lz$ estimates}\end{center}
In this section, we always assume that  the auxiliary $\rho(x)$ satisfies Lemma 2.1.
\begin{thm}\label{t3.1.}\hspace{-0.1cm}{\rm\bf 3.1.}\quad
Fix $1<q\le\fz, a\ge 1$, $\tz>0$ and $\wz\in RH_{s'}^{loc}$ with $1\le s<\fz$ and $1/s+1/s'=1$. Then, there exists $C=C(q,n,a,\wz,s,\tz)$
and $K_0=K_0(n,a)\ge 1$ with the following properties: Assume that $F, G, H_1$ are non-negative measurable functions on $\rz$ such that for any cube $Q=Q(x_0,r)$ with $r<\rho(x_0)$, there exist non-negative functions $G_{8Q}$ and $H_{8Q}$ with $F(x)\le G_{8Q}(x)+H_{8Q}(x)$ for ${\rm a.e.}\ x\in 8Q$ and
$$\l(\dfrac 1{|8Q|}\dint_{8Q} H_{8Q}(y)^qdy\r)^{1/q}\le a(M_{\rho,\tz}F(x)+H_1(\bar x)),\quad \forall x,\bar x\in 8Q;\eqno(3.1)$$
and for any $ x\in \rz$
$$\dsup_{x\in Q,r<\rho(x_0)}\dfrac
1{|8Q|}\dint_{8Q(x_0,r)}G_{8Q}(y)\,dy+ \dsup_{x\in Q,r\ge \rho(x_0)}\dfrac
1{\Psi_\tz(Q)|Q|}\dint_{Q(x_0,r)}|F(y)|\,dy\le G(x).\eqno(3.2)$$
Then for all $\lz>0$, $K\ge K_0$ and $0<\gz<1$
$$\wz\{M_{\rho,\tz} F>K\lz, G+H\le \gz\lz\}\le C\l(\dfrac {a^q}{K^q}+\dfrac \gz K\r)^{1/s}\wz\{M_{\rho,\tz} F>\lz\}.\eqno(3.3)$$
As a consequence, for all $0<p<q/s$, we have
$$\|M_{\rho,\tz} F\|_{L^p(\wz)}\le C(\|G\|_{L^p(\wz)}+\|H_1\|_{{L^p(\wz)}}),\eqno(3.4)$$
provided $\|M_{\rho,\tz} F\|_{L^p(\wz)}<\fz$, and
$$\|M_{\rho,\tz} F\|_{L^{p,\fz}(\wz)}\le C(\|G\|_{L^{p,\fz}(\wz)}+\|H_1\|_{{L^{p,\fz}(\wz)}}),\eqno(3.5)$$
provided $\|M_{\rho,\tz} F\|_{L^{p,\fz}(\wz)}<\fz$. Furthermore, if $p\ge 1$ then (3.4) and (3.5) hold, provided $F\in L^1$(whether or not $M_{\rho,\tz} F\in L^p(\wz)$).
\end{thm}
{\it Proof.}\quad We  borrow some ideas from \cite{am, t1}. It suffices to consider the case $H=G$: indeed, set $\wt G=G+H$. Then (3.1) holds with $\wt G$ in place of $H$ and also (3.2) holds with $\wt G$ in place of $G$.

Since from on we assume that $H=G$. Set $E_\lz=\{x\in\rz:\ M_{\rho,\tz}F(x)>\lz\}$ which is assumed to have finite measure
(otherwise there is nothing to prove). Clearly, $E_\lz$ is an open set. By Whitney's decomposition, there exists a family of pairwise disjoint cubes $\{Q_j\}_j$ so that $E_\lz=\bigcup_j Q_j$ with the property  that
$4Q_j$ meets $E^c_\lz$, that is , there exists $x_j\in 4Q_j$ such that
$M_{\rho,\tz}F(x_j)\le\lz.$

Set $B_\lz=\{M_{\rho,\tz}F>K\lz, 2G\le \gz \lz\}$. Since $K\ge 1$ we have that $B_\lz\subset E_\lz$. Therefore
$B_\lz\subset \bigcup_j B_\lz\bigcap Q_j$. We first claim that if $Q_j=Q(x^0_j,r_j)$ with $r_j\ge \rho(x^0_j)$, then
 $$B_\lz\bigcap Q_j=\emptyset.$$
 In fact, if $r\ge \rho(x_j)$ and $x\in B_\lz\bigcap Q_j$, then by (3.2),  we have
 $$K_0\lz<\dfrac
1{\Psi_\tz(Q_j)|Q_j|}\dint_{Q_j}|F(y)|\,dy\le G(x)<\gz\lz$$
 but $\gz<1$, hence, the set in question is  empty.

 Hence, we next only consider these cubes $Q_j=Q(x^0_j,r_j)$ with $r_j<\rho(x^0_j)$.
 For each $j$ we assume that $B_\lz\bigcap Q_j\not=\emptyset$(otherwise
we discard this cube) and so there is $\bar x_j\in Q_j$ so that $G(\bar x_j)\le \gz\lz/2$.
Since $M_{\rho,\tz}F(x_j)\le\lz$, there is $\wt C_0$ depending only on $n,\tz,l_0,C_0$ such that for every $K\ge \wt C_0$ we have
$$\begin{array}{cl}
|B_\lz\bigcap Q_j|&\le |\{M_{\rho,\tz}F>K\lz\}\bigcap Q_j|\le |\{M_{\rho,\tz}(F\chi_{8Q_j})>(K/\wt C_0)\lz\}|\\
& \le |\{M_{\rho,\tz}(G_{8Q_j}\chi_{8Q_j})>(2K/\wt C_0)\lz\}|+|\{M_{\rho,\tz}(H_{8Q_j}\chi_{8Q_j})>(2K/\wt C_0)\lz\}|,
\end{array}$$
where we have used $F\chi_{8Q_j}\le G_{8Q_j}\chi_{8Q_j}+H_{8Q_j}\chi_{8Q_j}$ a.e. and $\chi_{8Q_j}$ is the characteristic function of $8Q_j$.  Let $c_p$ be the weak-type $(p,p)$ bound of the maximal function $M_{\rho,\tz}$. From (3.2) and $\bar x_j\in Q_j\subset 8Q_j$, we obtain
$$\begin{array}{cl}
|\{M_{\rho,\tz}(G\chi_{8Q_j})>(2K/\wt C_0)\lz\}|&\le \dfrac{2\wt C_0c_1}{K\lz}\dint_{8Q_j}G_{8Q_j}(y)dy\\
&\le \dfrac{2C_0c_1}{K\lz}|8Q_j|G(\bar x_j)|\le \dfrac{8^n\wt C_0c_1}{K\lz}|Q_j|\gz.
\end{array}$$
Next assume first that $q<\fz$. By (3.1) and $x_j,\bar x_j\in 8Q_j$, we obtain
$$\begin{array}{cl}
|\{M_{\rho,\tz}(H\chi_{8Q_j})>(2K/\wt C_0)\lz\}|&\le \l(\dfrac {2\wt C_0c_q}{K\lz}\r)^q\dint_{8Q_j}H_{8Q_j}^q(y)dy\\
&\le \l(\dfrac {2\wt C_0c_q}{K\lz}\r)^q|8Q_j|a^q(M_{\rho,\tz}F(x_j)+G(\bar x_j))^q\\
&\le \l(\dfrac {4\wt C_0c_qa8^n}{K}\r)^q|Q_j|.
\end{array}$$
From the two inequalities above, we get
$$|B_\lz\bigcap Q_j|\le C\l(\dfrac {a^q}{K^q}+\dfrac\gz K\r)|Q_j|.$$
Note that $\wz\in RH_{s'}^{loc}$. If $s'<\fz$, for any cube $Q=Q(x_0,r)$ with $r<\rho(x_0)$ and any measurable set
 $E\subset Q$ we have
 $$\dfrac {\wz(E)}{\wz(Q)}\le \dfrac {|Q|}{\wz(Q)}\l(\dfrac 1{|Q|}\dint_Qw(y)^{s'}dy\r)^{ 1/s'}\l(\dfrac {|E|}{|Q|}\r)^{1/s}\le C_\wz\l(\dfrac {|E|}{|Q|}\r)^{1/s}.$$
 Note that the same conclusion holds in the case $s'=\fz$. Applying this to $B_\lz\bigcap Q_j\subset Q_j$, we have
 $$\wz(B_\lz\bigcap Q_j)\le C_\wz C\l(\dfrac {a^q}{K^q}+\dfrac\gz K\r)^{1/s}\wz(Q_j).$$
 Since the Whitney cubes are disjoint we get
 $$\wz(B_\lz)\le \dsum_j\wz(B_\lz \bigcap Q_j)\le C\l(\dfrac {a^q}{K^q}+\dfrac\gz K\r)^{1/s'}\dsum_j\wz(Q_j)
 =C\l(\dfrac {a^q}{K^q}+\dfrac\gz K\r)^{1/s}\wz(E_\lz)$$
 which is (3.3).

 When $q=\fz$, then by (3.1)
 $$\|M_{\rho,\tz}(H_{8Q_j\chi_{8Q_j}})\|_{L^\fz}\le \|H_{8Q_j\chi_{8Q_j}}\|_{L^\fz}\le a(M_{\rho,\tz}F(x_j)+G(\bar x_j))\le 2a\lz.$$
Thus, choosing $K\ge 4a\wt C_0$ it follows that $\{M_{\rho,\tz}(H_{8Q_j\chi_{8Q_j}})>(K/2\wt C_0)\lz\}=\emptyset$. Hence,
we can obtain the desired result (with $K^{-q}=0$).

When $M_{\rho,\tz}F\in L^p(\wz)$, we show (3.4). If $q<\fz$, integrating the two-parameter good-$\lz$ inequality (3.3) against $p\lz^{p-1}d\lz$ on $(0,\fz)$, for $0<p<\fz$,
$$\|M_{\rho,\tz}F\|^p_{L^p(\wz)}\le CK^p\l(\dfrac {a^q}{K^q}+\dfrac\gz K\r)^{1/s}\|M_{\rho,\tz}F\|^p_{L^p(\wz)}+\dfrac{2^pK^p}{\gz^p}\|G\|^p_{L^p(\wz)}.$$
For $0<p<q/s$ we can choose $K$ large enough and then $\gz$ small enough, we can get (3.4). In the same way, if $M_{\rho,\tz}F\in L^{p,\fz}(\wz)$, one shows the corresponding estimate in $L^{p,\fz}(\wz)$.

Observe that in  the case $q=\fz$, $K$ is already chosen and we only have to take some small $\gz$. Thus, the corresponding estimates holds for $0<p<\fz$ no matter the value of $s$.

Finally, we consider the case $p\ge 1$ and $F\in L^1$. By  the standard method in pages 247-248 in \cite{am}, we can obtain the desired result. \hfill$\Box$
\medskip

{\bf Remark 3.1.}\quad  In Theorem 3.1, if $q=\fz$, in fact, we only need $\wz\in A_{\rho,\fz}^{\fz}=\bigcup_{p\ge 1}A_p^{\rho,\fz}$, no matter the value of $s$. If $s>1$ and $q<\fz$, then one also obtains the end-point $p=q/s$.  In addition, it should be worth pointing out that Theorem 3.1 generalizes Theorem 2.1 in \cite{t1}.

Next we give a application of Theorem 3.1 toward weighted norm inequalities for operators, avoiding all use of kernel representation.
\begin{thm}\label{t3.2.}\hspace{-0.1cm}{\rm\bf 3.2.}\quad
Let $1\le p_0<q_0\le\fz$. Let $T$ is sublinear operator acting on  $L^{p_0}$. Let $\{A_r\}_{r>0}$ be family of operators acting from $L_c^\fz$ into $L^{p_0}$. Assume that for any $\tz>0$
$$\begin{array}{cl}
&\dsup_{x\in B,r<\rho(x_0)}\dfrac
1{|c_0B|}\dint_{c_0B(x_0,r)}|T(I-A_{r(c_0B)})f(y)|^{p_0}\,dy\\
&+ \dsup_{x\in B,r\ge \rho(x_0)}\dfrac
1{\Psi_\tz(B)|B|}\dint_{B(x_0,r)}|Tf(y)|^{p_0}\,dy\le C_\tz M_{\rho,\tz}(|f|^{p_0})(x),\ \forall x\in\rz
\end{array}\eqno(3.5)$$
and
$$\l(\dfrac 1{|c_0B|}\dint_{c_0B}T(A_{r(c_0B)})f(y)|^{q_0}dy\r)^{1/q_0} \le C_\tz M_{\rho,\tz}(|f|^{p_0})^{1/p_0}(x),\eqno(3.6)$$
for any ball $B=B(x_0,r)$ with $r<\rho(x_0)$ and all $x\in c_0B$ and $r(c_0B)$ denotes  $c_0B$ radius and $c_0$ is a constant depending only on $n$. Let $p_0<p<q_0$( or $p=q_0$ when $q_0<\fz$) and $\wz\in A^{\rho,\fz}_{p/p_0}\bigcap RH^{loc}_{(q_0/p)'}$. There  exists a constant $C$ such that
$$\|Tf\|_{L^p(\wz)}\le C\|Tf\|_{L^p(\wz)}\eqno(3.7)$$ for all $f\in L_c^\fz(\rz)$.
\end{thm}
{\it Proof.}\quad  We first notice that Theorem 3.1 still holds if the cubes $Q$ and  $8Q$ in the conditions 3.1  and 3.2 are replaced by the balls $B$ and $c_0B$ respectively, where $c_0$ is a constant depending only on $n$.  We now consider $q_0<\fz$ and $p_0<p\le q_0$. Let $f\in L_c^\fz$ and so $F=|Tf|^{p_0}\in L^1$. Fix a ball $B=B(x_0,r)$ with $r<\rho(x_0)$ . As $T$ is sublinear, we have
$$F\le G_{c_0B}+H_{c_0B}\equiv 2^{p_0-1}|T(I-A_{r(c_0B)}f|^{p_0}+2^{p_0-1}|TA_{r(c_0B)}f|^{p_0}.$$
Then (3.5) and (3.6) yield the corresponding conditions (3.1) and (3.2) with $q=q_0/p_0,$ $ H_1\equiv 0, a=2^{p_0-1}C^{p_0}$ and $G=2^{p_0-1}M_{\rho,\tz}(|f|^{p_0})$. Since $\wz\in RH^{loc}_{(q_0/p)'}$, then one pick $1<s<q_0/p$ such that $\wz\in RH^{loc}_{s'}$. Thus, Lemma 2.3 and Theorem 3.1  with $p/p_0>1$ in place of $p$ and $s=q_0/p$ yield
$$\begin{array}{cl}
\|Tf\|_{L^p(\wz)}^{p_0}&\le \|M_{\rho,\tz}F\|_{L^{p/p_0}(\wz)}\le C\|G\|_{L^{p/p_0}(\wz)}\\
&=C\|M_{\rho,\tz}(|f|^{p_0})\|_{L^{p/p_0}(\wz)}\le C\|f\|_{L^{p/p_0}(\wz)},\end{array}$$
where in the last estimate we have used  the fact that there exists $\tz_1>0$ such that $\wz\in A^{\rho,\tz_1}_{p/p_0}$ (since $\wz\in A^{\rho, \fz}_{p/p_0}$) and $\tz=\tz_1(p/p_0)'$.

In the case $q_0=\fz$ and $p<\fz$, Theorem 3.1 applies as before when $\wz\in A_{p/p_0}^{\rho,\fz}$ by Remark 3.1.\hfill$\Box$

A slight strengthening of the hypotheses in Theorem 3.1 furnishes weighted $L^p$ estimates for commutators $BMO_\tz(\rho)$ functions.
For any $k\in\nn$ we define the $k$th order commutator
$$T_b^kf(x)=T((b(x)-b)^kf)(x),\quad f\in L_c^\fz,\ x\in\rz.$$
Note that $T_b^0=T$. Commutators are usually considered for linear operators $T$ in which case they can be alternatively defined by recurrence: the first order commutators is
$$T_b^1f(x)=[b,T]f(x)=b(x)Tf(x)-T(bf)(x)$$ and for $k\ge 2$, the $k$th order commutators is given by
$T_b^k=[b, T_b^{k-1}]$.

We claim that since $T$ is bounded in $L^{p_0}$ then $T_b^kf$ is well defined in  $L_{loc}^q$ for any $0<q<p_0$ and for any $f\in L_c^\fz$: take a ball $B$ containing the support of $f$ and observe that by sublinearity for  a.e. $x\in\rz$
$$|T_b^kf(x)|\le \dsum_{m=0}^kC_{m,k}||b-b_B|^{k-m}|T((b-b_B)^mf)(x)|.$$
Lemma 2.4 implies
$$\dint_B|b(y)-b_Q|^{mp_0}|f(y)|^{p_0}dy\le C\|f\|_{L^\fz}^{p_0}\|b\|_{BMO_\tz}^{mp_0}|B|\Psi_{mp_0\tz'}(B)<\fz.$$
Hence, $T((b-b_B)^mf)\in L^{p_0}$ and the claim follows.
\begin{thm}\label{t3.3.}\hspace{-0.1cm}{\rm\bf 3.3.}\quad
Let $1\le p_0<q_0\le\fz$. Let $T$ is sublinear operator acting on  $L^{p_0}$. Let $\{A_r\}_{r>0}$ be family of operators acting from $L_c^\fz$ into $L^{p_0}$. Assume that for any $\tz>0$
$$\begin{array}{cl}
&\l(\dfrac
1{|c_0B|}\dint_{c_0B(x_0,r)}|T(I-A_{r(c_0B)})f(y)|^{p_0}\,dy\r)^{1/p_0}\\
&\qquad\qquad\le C_\tz\dsum_{j=1}^\fz\az_j
\l(\dfrac
1{\Psi_\tz(B_{j+1})|B_{j+1}|}\dint_{B_{j+1}}|f(y)|^{p_0}\,dy\r)^{1/p_0}\end{array}\eqno(3.7)$$ holds for any ball $B=B(x_0,r)$ with $r<\rho(x_0)$ and $B_{j+1}=2^{j+1}c_0B$,
$$\l(\dfrac
1{\Psi_\tz(B)|B|}\dint_{B(x_0,r)}|Tf(y)|^{p_0}\,dy\r)^{\frac 1{p_0}}\le C_\tz
\dsum_{j=1}^\fz\az_j\l(\dfrac
1{\Psi_\tz(B_{j+1})|B_{j+1}|}\dint_{B_{j+1}}|f(y)|^{p_0}\,dy\r)^{\frac 1{p_0}}\eqno(3.8)$$
holds for any ball $B=B(x_0,r)$ with $r\ge\rho(x_0)$ and $B_{j+1}=2^{j+1}B$,and
$$\l(\dfrac 1{|c_0B|}\dint_{c_0B}T(A_{r(c_0B)})f(y)|^{q_0}dy\r)^{\frac 1{q_0}} \le C_\tz \dsum_{j=1}^\fz\az_j
\l(\dfrac
1{\Psi_\tz(B_{j+1})|B_{j+1}|}\dint_{B_{j+1}}|Tf(y)|^{p_0}\,dy\r)^{\frac 1{p_0}},\eqno(3.9)$$
holds for any ball $B=B(x_0,r)$ with $r<\rho(x_0)$ and $B_{j+1}=2^{j+1}c_0B$. Let $p_0<p<q_0$( or $p=q_0$ when $q_0<\fz$) and $\wz\in A^{\rho,\fz}_{p/p_0}\bigcap RH^{loc}_{(q_0/p)'}$. If $\sum_j \az j^k<\fz$, then there  exists a constant $C$ such that for all $f\in L_c^\fz(\rz)$ and $b\in BMO_{\tz_1}(\rho)$,
$$\|T^k_bf\|_{L^p(\wz)}\le C\|b\|^k_{BMO_{\tz_1}(\rho)}\|Tf\|_{L^p(\wz)}.\eqno(3.10)$$
\end{thm}
{\it Proof.}\quad  We only prove the case $k=1$, the case $k\ge 2$ can be deduced by induction. Let us fix $p_0<p<q_0$ and $\wz\in A^{\rho,\tz_2}_{p/p_0}\bigcap RH^{loc}_{(q_0/p)'}$. We assume that $q_0<\fz$, for $q_0=\fz$ is similar.
Without loss of generality, $b\in BMO_{\tz_1}(\rho)\bigcap L^\fz$ and  $f\in L_c^\fz$, then  $|T_b^1f|^{p_0}\in L^1$. Set $F=|T_b^1f|^{p_0}$. In the proof, we always assume $\eta>(l_0+1)(\tz_1+\tz_2)$ large enough and $\eta_1=p_02^{[p_0]+1}\eta $.

  Given a ball $B=B(x_0,r)$, we first consider the case $r<
 \rho(x_0)$, we set $f_{\bar B,b}=(b_{4\bar B}-b)f$, $\bar B=c_0B$ and decompose $T^1_b$ as follows:
 $$\begin{array}{cl}
 |T_b^1f(x)|&=|T((b(x)-b)f)(x)|\le |b(x)-b_{\bar B}||Tf(x)|+|T((b_{4\bar B}-b)f)(x)|\\
 &\le |b(x)-b_{\bar B}||Tf(x)|+|T(I-A_{r(\bar B)})f_{B,b})(x)|+|TA_{r(\bar B)}f_{\bar B,b}(x)|.
 \end{array}$$
We observe that $F\le G_{\bar B}+H_{\bar B}$ where
$$G_{\bar B}=4^{p_0}(G_{\bar B,1}+G_{\bar B,2})=4^{p_0-1}(|b-b_{4\bar B}|^{p_0}|Tf|^{p_0}+|T(I-A_{r(\bar B)})f_{\bar B,b}|^{p_0})$$
and $H_{\bar B}=2^{p_0-1}|TA_{r(\bar B)}f_{\bar B,b}|^{p_0}$.
Fix any $x\in \bar B$, by (2.1), we have
 $$\dfrac 1{|\bar B|}\dint_{\bar B}G_{\bar B,1}(y)dy=\dfrac 1{|\bar B|}\dint_{\bar B}|b-b_{4\bar B}|^{p_0}|Tf|^{p_0}(y)dy
 \le C\|b\|_{BMO_{\tz_1}(\rho)}M_{\rho,\eta}^{[p_0]+2}(|Tf|^{p_0})(x).$$
  Using (3.7), (2.1),  and (2.3), and note that $\dsum_i\az_j j<\fz,$ we get
$$\begin{array}{cl}
\l(\dfrac 1{|\bar B|}\dint_{\bar B}G_{\bar B,2}(y)dy\r)^{1/p_0}&=\l(\dfrac 1{|\bar B|}\dint_{\bar B}T(I-A_{r(\bar B)})f_{\bar B,b}|^{p_0}(y)dy\r)^{1/p_0}\\
&\le C\dsum_{j=1}^\fz\az_j\l(\dfrac
1{\Psi_{\eta_1}(B_{j+1})|B_{j+1}|}\dint_{B_{j+1}}|f_{\bar B,b}|^{p_0}(y)dy\r)^{1/p_0}\\
&\le C \dsum_{j=1}^\fz\az_j\|(b-b_{4{\bar B}})/\Psi_{\eta_1}(B_{j+1})\|_{\exp L, B_{j+1}}M_{\rho,\eta}^{[p_0]+2}(|f|^{p_0})^{\frac 1{p_0}}(x)\\
&\le C \|b\|_{BMO_{\tz_1}(\rho)}M_{\rho,\eta}^{[p_0]+2}(|f|^{p_0})^{\frac 1{p_0}}(x)\dsum_{j=1}^\fz\az_j j\\
&\le C \|b\|_{BMO_{\tz_1}(\rho)}M_{\rho,\eta}^{[p_0]+2}(|f|^{p_0})^{\frac 1{p_0}}(x).
\end{array}$$
Hence,  we have
$$\begin{array}{cl}
\l(\dsup_{x\in\bar B,r<\rho(x_0)}\dfrac 1{|\bar B|}\dint_{\bar B}G_{\bar B}(y)dy\r)^{1/p_0}&\le C(M_{\rho,\eta}^{[p_0]+2}(|Tf|^{p_0})^{\frac 1{p_0}}(x)+M_{\rho,\eta}^{[p_0]+2}(|f|^{p_0})^{\frac 1{p_0}}(x))\\
&\equiv G(x).\end{array}$$
We next estimate the average of $H_{\bar B}^q$ on $\bar B$ with $q=q_0/p_0$. Using (3.9) with $\tz=\eta_1$, we get
$$\begin{array}{cl}
\l(\dfrac 1{|\bar B|}\dint_{\bar B}H_{\bar B}^q(y)dy\r)^{1/q_0}&\le C\l(\dfrac 1{|\bar B|}\dint_{\bar B}TA_{r(\bar B)}f_{\bar B,b}|^{p_0}(y)dy\r)^{1/p_0}\\
&\le C\dsum_{j=1}^\fz\az_j\l(\dfrac
1{\Psi_{\eta_1}(B_{j+1})|B_{j+1}|}\dint_{B_{j+1}}|Tf_{\bar B,b}|^{p_0}(y)dy\r)^{1/p_0}\\
&\le C(M_{\rho,\eta_1}F)^{\frac 1{p_0}}(x)\\
&\quad+C\dsum_{j=1}^\fz\az_j\l(\dfrac
1{\Psi_{\eta_1}(B_{j+1})|B_{j+1}|}\dint_{B_{j+1}}|b-b_{4\bar B}|^{p_0}|Tf|^{p_0}(y)dy\r)^{1/p_0}\\
&\le C(M_{\rho,\eta_1}F)^{\frac 1{p_0}}(x)\\
&\quad +C \dsum_{j=1}^\fz\az_j\|(b-b_{4{\bar B}})/\Psi_{\eta}(B_{j+1})\|_{\exp L, B_{j+1}}M_{\rho,\eta}^{[p_0]+2}(|f|^{p_0})^{\frac 1{p_0}}(\bar x)\\
&\le C(M_{\rho,\eta_1}F)^{\frac 1{p_0}}(x)+C\|b\|_{BMO_{\tz_1}(\rho)}M_{\rho,\eta}^{[p_0]+2}(|f|^{p_0})^{\frac 1{p_0}}(\bar x).
\end{array}$$
for any $x,\ \bar x\in \bar B$. Thus
$$\begin{array}{cl}
\l(\dfrac 1{|\bar B|}\dint_{\bar B}H_{\bar B}^q(y)dy\r)^{1/q_0}&\le C(M_{\rho,\eta_1}F)^{\frac 1{p_0}}(x)+M_{\rho,\eta}^{[p_0]+2}(|f|^{p_0})^{\frac 1{p_0}}(\bar x))\\
&\equiv C((M_{\rho,\eta_1}F)^{\frac 1{p_0}}(x)+H_1(\bar x)).
\end{array}$$
Given a ball $B=B(x_0,r)$, we first consider the case $r\ge
 \rho(x_0)$, we set $f_{ B,b}=(b_{4 B}-b)f$, $\bar B=B$ and decompose $T^1_b$ as follows:
 $$ |T_b^1f(x)|=|T((b(x)-b)f)(x)|\le |b(x)-b_{ B}||Tf(x)|+|T(f_{ B,b})(x)|.$$
Fix any $x\in  B$, by (2.1), we have
 $$\dfrac 1{| B|}\dint_{ B}|b-b_{4B}|^{p_0}|Tf|^{p_0}(y)dy
 \le C\|b\|_{BMO_{\tz_1}(\rho)}M_{\rho,\eta}^{[p_0]+2}(|Tf|^{p_0})(x).$$
 Using (3.8) with $\tz=\eta_1$, by (2.1) and (2.3),  we get
$$\begin{array}{cl}
&\l(\dfrac 1{\Psi_{\eta_1}(B)| B|}\dint_{ B}|T(f_{ B,b})|^{p_0}(y)dy\r)^{1/p_0}\\
&\qquad\le C\dsum_{j=1}^\fz\az_jj\l(\dfrac
1{\Psi_{\eta_1}(B_{j+1})|B_{j+1}|}\dint_{B_{j+1}}|f_{ B,b}|^{p_0}(y)dy\r)^{1/p_0}\\
&\qquad\le C \dsum_{j=1}^\fz \az_jj\|(b-b_{4{ B}})/\Psi_{\eta(B_{j+1})}\|_{\exp L, B_{j+1}}M_{\rho,\eta}^{[p_0]+2}(|f|^{p_0})^{\frac 1{p_0}}(x)\\
&\qquad\le C \|b\|_{BMO_{\tz_1}}M_{\rho,\eta}^{[p_0]+2}(|f|^{p_0})^{\frac 1{p_0}}(x)\dsum_{j=1}^\fz \az_jj\\
&\qquad\le C \|b\|_{BMO_{\tz_1}}M_{\rho,\eta}^{[p_0]+2}(|f|^{p_0})^{\frac 1{p_0}}(x),
\end{array}$$
 Hence,  we have
$$\begin{array}{cl}
\l(\dsup_{x\in B,r\ge \rho(x_0)}\dfrac 1{\Psi_{\eta_1}(B)| B|}\dint_{ B}F(y)dy\r)^{\frac 1{p_0}}&\le C(M_{\rho,\eta}^{[p_0]+2}(|Tf|^{p_0})^{\frac 1{p_0}}(x)+M_{\rho,\eta}^{[p_0]+2}(|f|^{p_0})^{\frac 1{p_0}}(x))\\
&\equiv G(x).
\end{array}$$
Thus,
applying Lemma 2.3 and Theorem 3.1, we get
$$\begin{array}{cl}
\|T_b^1f\|_{L^p(\wz)}^{p_0}&\le \|M_{\rho,\eta_1}F\|_{L^{p/p_0}(\wz)}
\le C\|G\|_{L^{p/p_0}(\wz)}+C\|H_1\|_{L^{p/p_0}(\wz)}\\
&\le C\|M_{\rho,\eta}^{[p_0]+2}(|f|^{p_0})\|_{L^{p/p_0}(\wz)}+C\|M_{\rho,\eta}^{[p_0]+2}(|Tf|^{p_0})\|_{L^{p/p_0}(\wz)}\\
&\le C\|f\|_{L^p(\wz)}^{p_0}+\|Tf\|_{L^p(\wz)}^{p_0}\le C\|f\|_{L^p(\wz)}^{p_0},
\end{array}$$
if  $\eta$ is large enough.\hfill$\Box$

\begin{center} {\bf 4. Spectral multipliers  }\end{center}
Suppose that $L$ is a nonnegative self-adjoint operator acting on $L^2(\rz)$. Let $E(\lz)$ be the spectral resolution of $L$. By the spectral
theorem, for any bounded Borel function $F:  [0,\fz)\to \cc$ one can define the operator
$$F(L)=\dint_0^\fz F(\lz)dE(\lz),$$
which is bounded on $L^2(\rz)$. The question of $L^p$ estimates for functions of a self-adjoint operator is a delicate one. In fact, even for a Schr\"odinger operator $H=-\triangle+V(x)$ with a nonnegative potential, and a bound smooth kernel and hence does not fall within the scope of the Calder$\acute{o}$n-Zygmund theory. The first to overcome this difficulty was Hebish \cite{h}. Later, J. Dziuba\'{n}ski \cite{dz} gave
 a spectral multiplier theorem for $H^1$ spaces associated with Schr\"odinger operators with potentials satisfying a reverse H\"older inequality. On the other hand, X. T Doung, etc \cite{dos} showed that a sharp spectral multiplier for a non-negative self-adjoint operator $L$ was obtained under the assumption of the kernel $p_t(x,y)$ of the analytic semigroup $e^{-tL}$ having a Gaussian upper bound. Recently, T Doung, etc \cite{dsy} generalized the main results \cite{dos} to weighted cases; see also \cite{cd}.

A natural problem  considered in the spectral multiplier theory is to
give sufficient conditions on $F$ and $L$ which imply the weighted boundedness of $F(L)$  associated with Schr\"odinger operators.

In this section, we always assume that $L$ is a non-negative self-adjoint operator on  $L^2(\rz)$ and that the  semigroup $e^{-tL}$, generated by $-L$ on $L^2(\rz)$, has the kernel $p_t(x,y)$  which satisfies the following Gaussian upper bound
$$|p_t(x,y)|\le C_N t^{-n/2}\l(1+\dfrac{\sqrt{t}}{\rho(x)} +\dfrac{\sqrt{t}}{\rho(y)}\r)^{-N}\exp{\l(-b\dfrac{|x-y|^2}{t}\r)}\eqno(4.1)  $$
for all $t>0$ , $N>0$, and $x,y\in\rz$, where the auxiliary function $\rho(x)$ satisfies Lemma 2.1, $C_N$ depends only $N$, and $b$ is a positive constant.

Such estimates are typical for divergence Schr\"odinger differential operator
$$L=-\pz_i(a_{ij}(x)\pz_j)+V(x)\ {\rm on}\ \rz,\ \ n\ge 3,$$ where  $V(x)\in RH_{n/2}$
is a  nonnegative potential, and $a_{ij}$ satisfy (1.1) (see \cite{dz}).

Suppose that $T$ is a bounded operator on $L^2$. We say that a measurable function $K_T: \rr^{2n}\to\cc$ is the (singular) kernel of $T$ if
$$<Tf_1, f_2>=\dint_\rz Tf_1\bar{f_2}dx=\dint_\rz K_T(x,y)f_1(y)\overline{f_2(x)}dxdy$$ for all $f_1,f_2\in C_c(\rz)$( for all $f_1,f_2\in C_c(\rz)$ such that $\supp f_1\bigcap\supp f_2=\emptyset$, respectively).

\begin{thm}\label{t4.1.}\hspace{-0.1cm}{\rm\bf 4.1.}\quad
Let $T$ be a non-negative self-adjoint operator such that the corresponding heat kernel satisfy (4.1). Suppose
that $F [0,\fz)\to\cc$ is a bounded Borel function such that
$$\dsup_{t>0}\|\eta\dz_tF\|_{W_s^\fz}\le C_s<\fz,\eqno(4.2)$$
for any $s>0$, where $\dz_tF(\lz)=F(t\lz),\ \|F\|_{W_s^\fz}=\|(I-d^2/d^2)^{s/2}F\|_{L^\fz}$ and $\eta\in C_c^\fz(\rr_+)$ is a fixed
function, non identically zero. Then the operator $F(L)$ is bounded on $L^p(\wz)$ for all $p$ and $\wz$ satisfying $1<p<\fz$ and $\wz\in A_p^{\rho,\fz}$, and is of weighted  weak type (1,1) for weights $\wz\in A_1^{\rho,\fz}$. In addition,
$$\|F(L)f\|_{L^p(\wz)\to L^p(\wz)}\le C_s\l(\dsup_{t>0}\|\eta\dz_tF\|_{W_s^\fz}+|F(0)|\r),\eqno(4.3)$$
for $1<p<\fz$, $\wz\in A_p^{\rho,\fz}$ provided that $s$ is large enough, and
$$\|F(L)f\|_{L^p(\wz)\to L^p(\wz)}\le C_s\l(\dsup_{t>0}\|\eta\dz_tF\|_{W_s^\fz}+|F(0)|\r),\eqno(4.4)$$
for $\wz\in A_1^{\rho,\fz}\bigcap RH_2^{loc}$ provided that $s$ is large enough.
\end{thm}
For the commutators for $F(L)$, we have the following results.
\begin{thm}\label{t4.2.}\hspace{-0.1cm}{\rm\bf 4.2.}\quad
Let $T$ be a non-negative self-adjoint operator such that the corresponding heat kernel satisfy (4.1). Suppose
that $F [0,\fz)\to\cc$ is a bounded Borel function such that
$$\dsup_{t>0}\|\eta\dz_tF\|_{W_s^\fz}\le C_s<\fz,$$
for any $s>0$, where $\dz_tF(\lz)=F(t\lz),\ \|F\|_{W_s^\fz}=\|(I-d^2/d^2)^{s/2}F\|_{L^\fz}$ and $\eta\in C_c^\fz(\rr_+)$ is a fixed
function, non identically zero. Let $b\in BMO_\tz(\rho)$ and $k\in\nn$, then for $k$-order commutator $F_b^k(L)$, we
have
$$\|F_b^k(L)f\|_{L^p(\wz)}\le C \|b\|_{BMO_\tz(\rho)}^k \|f\|_{L^p(\wz)},\eqno(4.5)$$
for $1<p<\fz$, $\wz\in A_p^{\rho,\fz}$, and
$$\wz(\{x\in\rz: \ |F_b^k(L)f(x)|>\lz\})\le C\Phi(\|b\|_{BMO_\tz(\rho)}) \dint_\rz\Phi\l(\dfrac
{|f(x)|}\lz\r)\wz(x)dx,\eqno(4.6)$$ where
$\Phi(t)=t\log(e+t)^k$ and $\wz\in A_1^{\rho,\fz}\bigcap RH_2^{loc}$.
\end{thm}

{\bf Remark 4.1.}\quad Let $F(\lz)=\lz^{i\gz}$ with $\gz\in\rr$, $F(\lz)=\cos\lz$, $F(\lz)=\sin\lz$, or $F(\lz)\in C_0^\fz(0,\fz)$, then these  $F$ all satisfy (4.2).

{\bf Remark 4.2.}\quad In fact, Theorems 3.1, 3.2, 3.3 and Theorems 4.1 and 4.2 still hold on RD spaces; see \cite{yz} for more details.

\subsection*{4.1. Some lemmas}
To prove Theorems 4.1 and 4.2, we need some lemmas.

\begin{lem}\label{l4.1.}\hspace{-0.1cm}{\rm\bf 4.1.}\quad
Suppose that (4.1) holds. Then for any $N>0$ such that
$$\dint_{\rz\setminus B(y,r)}|p_t(x,y)|^2dx\le Ct^{-\frac n2}\l(1 +\dfrac{\sqrt{t}}{\rho(y)}\r)^{-N}\exp{\l(-b\dfrac{r^2}{t}\r)}  $$
\end{lem}
The proof is obviously by (4.1).

\begin{lem}\label{l4.2.}\hspace{-0.1cm}{\rm\bf 4.2.}\quad
Suppose that (4.1) holds. For any $s\ge 0$ and $N>0$ there exists a constant $C$ such that
$$\dint_\rz|p_{(1+i\tau)R^{-2}}(x,y)|^2|x-y|^sdx\le CR^{n-s}(1+|\tau|)^s\l(1 +\dfrac{R}{\rho^2(y)}\r)^{-N}.\eqno(4.7)$$
where $p_{(1+i\tau)R^{-2}}=K_{\exp((1+i\tau)R^{-2}L)}$.
\end{lem}
{\it Proof.}\quad Assume $\|f\|_{L^2}$ and $f\in \rz\setminus B(y,r)$. we defines the holomorphic function $F_y:\{z\in\cc:{\cal R}e\,z>0\}\to\cc$ by the formula
$$F_y(z)=e^{-zR^2}R^n \l(1 +\dfrac{z}{\rho^2(y)}\r)^{-N}\l(\dint_\rz p_z(x,y)f(x)dx\r)^2.$$
If let $z=|z|e^{i\tz}$, then
$$\|p_z(\cdot,y)\|_{L^2}^2=\|p_{|z|\cos\tz}(\cdot,y)\|_{L^2}^2.$$
From this and Lemma 4.1, we get
$$\begin{array}{cl}
|F_y(z)|&\le e^{-R^2|z|\cos\tz}R^{-n}\|p_{|z|\cos\tz}(\cdot,y)\|_{L^2}^2\\
&\le e^{-R^2|z|\cos\tz}R^{-n}(|z|\cos\tz)^{-\frac 12}\\
&\le CR^{-n}(|z|\cos\tz)^{-\frac n2}.
\end{array}$$
Similarly, for $\tz=0$ by Lemma 4.1,
$$|F_y(|z|)|\le CR^{_n}|z|^{-n/2}e^{-br^2/|z|}.$$
Combining the inequalities above and Lemma 9 in \cite{d}(see also the proof of Lemma 4.1 in \cite{dos}), we have
$$|F_y((1+i\tau)R^{-2})|\le Ce^{-b(rR/(1+|\tau|))^2}.$$
$$\dint_\rz|p_{(1+i\tau)R^{-2}}(x,y)|^2dx\le CR^n\l(1 +\dfrac1{R\rho(y)}\r)^{-N}e^{-b(rR/(1+|\tau|))^2}.$$
Hence, we have
$$\begin{array}{cl}
\dint_\rz&|p_{(1+i\tau)R^{-2}}(x,y)|^2|x-y|^sdx\\
&=\dsum_{k\ge 0}\dint_{k(1+|\tau|)R^{-1}\le|x-y|\le (k+1)(1+|\tau|)R^{-1}}|p_{(1+i\tau)R^{-2}}(x,y)|^2|x-y|^sdx\\
&\le (1+|\tau|)^sR^{-s}\dsum_{k\ge 0}(k+1)^s\dint_{\rz\setminus B(y,(k+1)(1+|\tau|)R^{-1})}|p_{(1+i\tau)R^{-2}}(x,y)|^2dx\\
&\le CR^{n-s}(1+|\tau|)^s\l(1 +\dfrac1{R\rho(y)}\r)^{-N}.
\end{array}$$
Thus, (4.7) holds.
\hfill$\Box$

Applying Lemmas 4.1 and 4.2, and adapting the similar arguments in the proof of lemma 4.2 in \cite{dos}, we can obtain the following result.
\begin{lem}\label{l4.3.}\hspace{-0.1cm}{\rm\bf 4.3.}\quad
Let $R>0,\ s>0$. Then for any $\ez>0$ and $N>0$, there exists a constant $C=C(s,\ez,N)$ such that
$$\dint_\rz |K_{F(\sqrt{L})}|^2(1+R|x-y|)^sdx\le CR^n\l(1 +\dfrac1{R\rho(y)}\r)^{-N}\|\dz_R F\|_{W^\fz_{s/2+\ez}}^2$$
for all Borel functions $F$ such that $\supp\ F\subset [R/4,R]$.
\end{lem}

\begin{lem}\label{l4.4.}\hspace{-0.1cm}{\rm\bf 4.4.}\quad
Suppose $\wz\in A^{\rho,\tz}_1\bigcap RH_2^{loc}$. Let $\eta=3\tz+(l_0+2)n$, then for any $s>n+2\eta$ such that
$$\dint_{\rz\setminus B(y,r)}(1+R|x-y|)^{-s}\wz(x)^2dx\le CR^{-n}(1+rR)^{n+2\eta-s}\l(1+\frac 1{r\rho(y)}\r)^{-2\eta}\wz(y)^2.$$
\end{lem}
{\it Proof.}\quad Let $\eta_1=2\tz+(l_0+1)n$. Assume that $rR>1$. Then
$$\begin{array}{cl}
\dint_{\rz\setminus B(y,r)}&(1+R|x-y|)^{-s}\wz(x)^2dx\\
&\le \dsum_{k\ge 0}\dint_{2^kr\le |x-y|\le 2^{k+1}r}
(R|x-y|)^{-s}\wz(x)^2dx\\
&\le \dsum_{k\ge 0}(2^krR)^{-s+n}\dfrac 1{(2^{k+1}r)^n}\dint_{ |x-y|\le 2^{k+1}r}
\wz(x)^2dx\\
&\le C\dsum_{k\ge 0}(2^krR)^{-s+n}\l(\dfrac 1{(2^{k+1}r)^n}\dint_{ |x-y|\le 2^{k+1}r}
\wz(x)dx\r)^2\l(1+\dfrac {2^{k+1}r}{\rho(y)}\r)^{2\eta_1}\\
&\le C\dsum_{k\ge 0}(2^krR)^{-s+n}(M_{\rho,\tz}w(y))^2\l(1+\dfrac {2^{k+1}r}{\rho(y)}\r)^{2\eta}\\
&\le C\dsum_{k\ge 0}(2^krR)^{-s+n}w(y)^2\l(1+\dfrac {2^{k+1}rR}{R\rho(y)}\r)^{2\eta}\\
&\le C\dsum_{k\ge 0}(2^krR)^{-s+n+2\eta}w(y)^2\l(1+\dfrac 1{R\rho(y)}\r)^{2\eta}\\
&\le C(rR)^{-s+n+2\eta}w(y)^2\l(1+\dfrac 1{R\rho(y)}\r)^{2\eta},
\end{array}$$
since $s>n+2\eta$.

If $rR<1$, note that $s>n+2\eta$, we then have
$$\begin{array}{cl}
\dint_{\rz}(1+R|x-y|)^{-s}\wz(x)^2dx&\le \dint_{|x-y|<1/R}\wz(x)^2dx\\
&+\dsum_{k\ge 1}2^{-ks}\dint_{2^k/R\le |x-y|\le 2^{k+1}/R}
\wz(x)^2dx\\
&\le C R^{-n}\dsum_{k\ge 0}2^{k(-s+n+2\eta)}w(y)^2\l(1+\dfrac 1{R\rho(y)}\r)^{2\eta}\\
&\le C R^{-n}w(y)^2\l(1+\dfrac 1{R\rho(y)}\r)^{2\eta}.
\end{array}$$\hfill$\Box$

Applying Lemmas 4.3 and 4.4, and adapting the similar arguments in the proof of Theorem 3.1 in \cite{dos}, we can obtain the following result.
\begin{lem}\label{l4.5.}\hspace{-0.1cm}{\rm\bf 4.5.}\quad
Suppose $\wz\in A^{\rho,\tz}_1\bigcap RH_2^{loc}$ and $\eta=3\tz+(l_0+2)n$. If for any $s>n/2+\eta$ such that  $\dsup_{t>0}\|\eta\dz_t F\|_{W^\fz_{s}}\le C_s$, then
$$\dint_{\rz\setminus B(y,r)}|K_{F(1-\Phi_r)(\sqrt{L})}(x,y)|\wz(x)dx\le C\wz(y),\quad {\rm a.e.}\ y\in\rz,\eqno(4.8)$$
where $\Phi_r(\lz)=\exp(-(\lz r)^2)$.
\end{lem}
\begin{lem}\label{l4.6.}\hspace{-0.1cm}{\rm\bf 4.6.}\quad
Suppose $\wz\in A^{\rho,\tz}_1\bigcap RH_2^{loc}$ and $b\in BMO_{\tz_1}(\rho)$. Let $\eta=(l_0+1)\tz_1+3\tz+(l_0+2)n$, then there exists $s>n+2\eta$ such that
$$\begin{array}{cl}
\dint_{\rz\setminus B(y,r)}&|b(x)-b_B|^2(1+R|x-y|)^{-s}\wz(x)^2dx\\
&\le C\|b\|^2_{BMO_{\tz_1}(\rho)}R^{-n}(1+rR)^{n+2\eta-s}\l(1+\dfrac 1{r\rho(y)}\r)^{-2\eta}\wz(y)^2,
\end{array}$$
where $y\in B=B(x_0,r)$ with $r< \rho(x_0)$ and $b_B$ is the average on $B$.
\end{lem}
{\it Proof.}\quad Let $\eta_1=2\tz+(l_0+1)(n+\tz_1)$. Assume that $rR>1$. Since $\wz\in  RH_2^{loc}$, So $\wz\in  RH_{2\gz}^{loc}$ for some $\gz>1$. Let $1/\gz+1/\gz'=1$. Then
$$\begin{array}{cl}
&\dint_{\rz\setminus B(y,r)}|b(x)-b_B|^2(1+R|x-y|)^{-s}\wz(x)^2dx\\
&\quad \le \dsum_{k\ge 0}\dint_{2^kr\le |x-y|\le 2^{k+1}r}
|b(x)-b_B|^2(R|x-y|)^{-s}\wz(x)^2dx\\
&\quad \le \dsum_{k\ge 0}(2^krR)^{-s+n}\dfrac 1{(2^{k+1}r)^n}\dint_{ |x-y|\le 2^{k+1}r}
|b(x)-b_B|^2\wz(x)^2dx\\
&\quad \le \dsum_{k\ge 0}(2^krR)^{-s+n}\l(\dfrac 1{(2^{k+1}r)^n}\dint_{ |x-y|\le 2^{k+1}r}
|b(x)-b_B|^{2\gz'}dx\r)^{1/\gz'}\\
&\qquad\qquad\times\l(\dfrac 1{(2^{k+1}r)^n}\dint_{ |x-y|\le 2^{k+1}r}
(\wz(x))^{2\gz}dx\r)^{1/\gz}\\
&\quad \le C\|b\|^2_{BMO_{\tz_1}(\rho)}\dsum_{k\ge 0}(2^krR)^{-s+n}\l(\dfrac 1{(2^{k+1}r)^n}\dint_{ |x-y|\le 2^{k+1}r}
\wz(x)dx\r)^2\l(1+\dfrac {2^{k+1}r}{\rho(y)}\r)^{2\eta_1}\\
&\quad \le C\|b\|^2_{BMO_{\tz_1}(\rho)}\dsum_{k\ge 0}(2^krR)^{-s+n}(M_{\rho,\tz}w(y))^2\l(1+\dfrac {2^{k+1}r}{\rho(y)}\r)^{2\eta}\\
&\quad \le C\|b\|^2_{BMO_{\tz_1}(\rho)}\dsum_{k\ge 0}(2^krR)^{-s+n}w(y)^2\l(1+\dfrac {2^{k+1}rR}{R\rho(y)}\r)^{2\eta}\\
&\quad \le C\|b\|^2_{BMO_{\tz_1}(\rho)}\dsum_{k\ge 0}(2^krR)^{-s+n+2\eta}w(y)^2\l(1+\dfrac 1{R\rho(y)}\r)^{2\eta}\\
&\quad \le C\|b\|^2_{BMO_{\tz_1}(\rho)}(rR)^{-s+n+2\eta}w(y)^2\l(1+\dfrac 1{R\rho(y)}\r)^{2\eta}.
\end{array}$$
If $rR<1$, similarly, we have
$$\begin{array}{cl}
\dint_{\rz}|b(x)-b_B|^2&(1+R|x-y|)^{-s}\wz(x)^2dx\\
&\le \dint_{|x-y|<1/R}|b(x)-b_B|^2\wz(x)^2dx\\
&\qquad+\dsum_{k\ge 1}2^{-ks}\dint_{2^k/R\le |x-y|\le 2^{k+1}/R}
|b(x)-b_B|^2\wz(x)^2dx\\
&\le C\|b\|^2_{BMO_{\tz_1}(\rho)} R^{-n}\dsum_{k\ge 0}2^{k(-s+n+2\eta)}w(y)^2\l(1+\dfrac 1{R\rho(y)}\r)^{2\eta}\\
&\le C\|b\|^2_{BMO_{\tz_1}(\rho)} R^{-n}w(y)^2\l(1+\dfrac 1{R\rho(y)}\r)^{2\eta}.
\end{array}$$\hfill$\Box$

Applying Lemmas 4.3 and 4.6, and adapting the similar arguments in the proof of Theorem 3.1 in \cite{dos}, we can obtain the following result.
\begin{lem}\label{l4.7.}\hspace{-0.1cm}{\rm\bf 4.7.}\quad
Suppose $\wz\in A^{\rho,\tz}_1\bigcap RH_2^{loc}$, $b\in BMO_{\tz_1}(\rho)$ and $\eta=(l_0+1)\tz_1+3\tz+(l_0+2)n$. If for any $s>n/2+\eta$ such that  $\dsup_{t>0}\|\eta\dz_t F\|_{W^\fz_{s}}\le C_s$, then
$$\dint_{\rz\setminus B(y,r)}|K_{F(1-\Phi_r)(\sqrt{L})}(x,y)||b(x)-b_B|\wz(x)dx\le C\|b\|_{BMO_{\tz_1}(\rho)}\wz(y),$$
where  $\Phi_r(\lz)=\exp(-(\lz r)^2)$, $y\in B=B(x_0,r)$ with $r\le \rho(x_0)$ and $b_B$ is the average on $B$.
\end{lem}

\begin{lem}\label{l4.8.}\hspace{-0.1cm}{\rm\bf 4.8(\cite{t1}).}\quad
For any a  ball $B=B(x_0,r)$, if $r\ge \rho(x_0)$, then the ball
$B$ can be decomposed into finite disjoint cubes $\{Q_i\}_{i=1,m}$
such that $B\subset \bigcup_i^{m} Q_i\subset 2\sqrt{n}B$ and
$r_i/2\le \rho(x)\le 2\sqrt{n}C_0r_i$ for some $x\in
Q_i=Q(x_i,r_i)$, where $C_0$ is  same as Lemma 2.1.
\end{lem}

\subsection*{4.2. Proof of Theorem 4.1}
In this section, we  borrow some ideas from \cite{dsy}.
 We first prove (4.5).   Since $\wz\in A_p^{\rho,\fz}$, then there exist $p_0>1$ and $\tz>0$ such that $\wz\in A_{p/p_0}^{\rho,\tz}$ by Proposition (iii).  We will show that for any $\eta>0$
$$\begin{array}{cl}
&\dsup_{x\in B,r<\rho(x_0)}\dfrac
1{|c_0B|}\dint_{c_0B(x_0,r)}|T(I-A_{r(c_0B)})f(y)|^{p_0}\,dy\\
&+ \dsup_{x\in B,r\ge \rho(x_0)}\dfrac
1{\Psi_\eta(B)|B|}\dint_{B(x_0,r)}|Tf(y)|^{p_0}\,dy\le C M_{\rho,\eta}(|f|^{p_0})(x),\ \forall x\in\rz
\end{array}\eqno(4.9)$$
for all $f\in L_c^\fz(\rz)$.

Let us now prove (4.9). Observe that  $\dsup_{t>0}\|\eta\dz_t F\|_{W^\fz_{s}}\sim  \dsup_{t>0}\|\eta\dz_t G\|_{W^\fz_{s}}$ where $G(\lz)=F(\sqrt{\lz})$. So, we can replace $F(L)$ by $F(\sqrt{L})$ in the proof. Notice that $F(\lz)=F(\lz)-F(0)+F(0)$ and hence
$$F(\sqrt{L})=(F(\cdot)-F(0))(\sqrt{L})+F(0)I.$$
Replacing $F$ by $F-F(0)$, we may assume in the sequel that $F(0)=0$. Let $\vz\in C_c^\fz(0,\fz)$ be a non-negative function satisfying $\supp\ \vz\subset [\frac 14,1]$ and $\sum_{l=-\fz}^\fz\vz(2^{-l}\lz)=1$ for any $\lz>0$, and let $\vz_l$ denote the function $\vz(2^{-l}\cdot)$. Then
$$F(\lz)=\dsum_{l=-\fz}^\fz\vz(2^{-l}\lz)F(\lz)=\dsum_{l=-\fz}^\fz F^l(\lz),\quad \forall \lz\ge 0.\eqno(4.10)$$
This decomposition  implies that the sequence $\sum_{l=-N}^NF^l(\sqrt{\lz})$ converges strongly in $L^2(\rz)$ to $F(\sqrt{\lz})$.

We first consider the case $B=B(x_0,r)$ with $r<\rho(x_0)$. For every $l\in\zz, r>0$, $M\in\nn$ and $\lz>0$, we set
$$F_{r,M}(\lz)=F(\lz)(1-e^{-(r\lz)^m})^M,\eqno(4.11)$$
$$F^l_{r,M}(\lz)=F(\lz)(1-e^{-(r\lz)^m})^M,\eqno(4.12)$$
We use the decomposition $f=\sum_{j=0}^\fz f_j$ in which $f_j=f\chi_{ U_j(\bar B)}$, where $U_0(B)=c_0B:=\bar B$ and $U_j(\bar B)=2^j\bar B\setminus 2^{j-1}\bar B$ for $j=1,2\cdots$. We set $r_B=c_0r$, then
$$\begin{array}{cl}
F(\sqrt{L})(1-e^{-r_B^2L})^Mf&=F_{r_B, M}(\sqrt{L})f\\
&=\dsum_{j=1}^2F_{r_B, M}(\sqrt{L})f_j+\dlim_{N\to \fz}\dsum_{l=-N}^N\dsum_{j=3}^\fz F_{r_B, M}^l(\sqrt{L})f_j,
\end{array}\eqno(4.13)$$
where the sequence converges strongly in $L^2(\rz)$.

Note that $\|e^{-tL}f\|_{L^p(\rz)}\le C\|f\|_{L^p(\rz)}$ for any $t>0$, and the $L^p$-boundedness of the operator $ F(\sqrt{L})$(see Theorem 3.1 in \cite{dos}), for any $x\in B$, we have
$$\begin{array}{cl}
\l(\dfrac
1{|c_0B|}\dint_{c_0B(x_0,r)}|F_{r_B, M}(\sqrt{L})f_j|^{p_0}\,dy\r)^{1/p_0}
&\le C|B|^{-1/p_0}\|F_{r_B, M}(\sqrt{L})f_j\|_{L^{p_0}(\rz)}\\
&\le C M_{\rho,\eta}(|f|^{p_0})^{1/p_0}(x),
\end{array}\eqno(4.14)$$
for $j=1,2$.

Fix $j\ge 3$. Let $p_1\ge 2$ and $\frac 1{p_0}-\frac 1{p_1}=\frac 12$. Adapting the same arguments in pages 1117-1119 in \cite{dsy}, using Lemma 4.3 with $N=0$,  we have
$$\begin{array}{cl}
&\l(\dfrac
1{|c_0B|}\dint_{c_0B(x_0,r)}|F^l_{r_B, M}(\sqrt{L})f_j|^{p_0}\,dy\r)^{1/p_0}\\
&\qquad\le C|B|^{-1/p_1}\|F^l_{r_B, M}(\sqrt{L})f_j\|_{L^{p_1}(\rz)}\\
&\qquad\le C|B|^{-1/p_1}\|F^l_{r_B, M}(\sqrt{L})\|_{L^{p_0}(U_j(\bar B))\to L^{p_1}(\bar B)}\|f_j\|_{L^{p_0}(\rz)}\\
&\qquad\le C 2^{j(\eta+n)/p_0}|B|^{\frac 12}\|F^l_{r_B, M}(\sqrt{L})\|_{L^{p_0}(U_j(\bar B))\to L^{p_1}(\bar B)}M_{\rho,\eta}(|f|^{p_0})^{1/p_0}(x)\\
&\qquad\le C 2^{-js+j(\eta+n)/p_0}\l(\min\{1,(2^lr_B)^{2M}\}\max\{1,(2^lr_B)^{\frac n2}\}\r)\\
&\qquad\qquad\times M_{\rho,\eta}(|f|^{p_0})^{1/p_0}(x)\dsup_{l\in\zz}\|\dz_{2^l}[\vz_lF]\|_{W^\fz_s}.
\end{array}\eqno(4.15)$$
Hence,
$$\begin{array}{cl}
&\dsum_{j=3}^\fz\dsum_{l=-\fz}^\fz\l(\dfrac
1{|c_0B|}\dint_{c_0B(x_0,r)}|F^l_{r_B, M}(\sqrt{L})f_j|^{p_0}\,dy\r)^{1/p_0}\\
&\qquad\le C \dsum_{j=3}^\fz 2^{-js+j(\eta+n)/p_0}\l(\dsum_{l=-\fz}^\fz\min\{1,(2^lr_B)^{2M}\}\max\{1,(2^lr_B)^{\frac n2}\}\r)\\
&\qquad\qquad\times M_{\rho,\eta}(|f|^{p_0})^{1/p_0}(x)\dsup_{l\in\zz}\|\dz_{2^l}[\vz_lF]\|_{W^\fz_s}\\
&\qquad\le C \dsum_{j=3}^\fz 2^{-js+j(\eta+n)/p_0}\l(\dsum_{l: 2^lr_B>1}(2^lr_B)^{-s+\frac n2}+\dsum_{l: 2^lr_B\le 1}(2^lr_B)^{2M-s}\r)\\
&\qquad\qquad\times M_{\rho,\eta}(|f|^{p_0})^{1/p_0}(x)\dsup_{l\in\zz}\|\dz_{2^l}[\vz_lF]\|_{W^\fz_s}\\
&\le CM_{\rho,\eta}(|f|^{p_0})^{1/p_0}(x)\dsup_{l\in\zz}\|\dz_{2^l}[\vz_lF]\|_{W^\fz_s},
\end{array}\eqno(4.15)$$
if $s>n+\eta$ and $M>s/2$.

We now consider the case $B=B(x_0,r)$ with $r\ge \rho(x_0)$. Let $f=\sum_{j=0}^\fz f_j$ in which $f_j=f\chi_{ U_j( B)}$ where $U_0(B)= B$ and $U_j( B)=2^jB\setminus 2^{j-1} B$ for $j=1,2\cdots$. We write
$$F(\sqrt{L})f
=\dsum_{j=1}^2F(\sqrt{L})f_j+\dlim_{N\to \fz}\dsum_{l=-N}^N\dsum_{j=3}^\fz F^l(\sqrt{L})f_j,
$$
where the sequence converges strongly in $L^2(\rz)$.

Similar to the proof of (4.14), we have
$$\begin{array}{cl}
\l(\dfrac
1{\Psi_\eta(B)|B|}\dint_{B(x_0,r)}|F(\sqrt{L})f_j|^{p_0}\,dy\r)^{1/p_0}
&\le C(\Psi_\eta(B)|B|)^{-1/p_0}\|F(\sqrt{L})f_j\|_{L^{p_0}(\rz)}\\
&\le C M_{\rho,\eta}(|f|^{p_0})^{1/p_0}(x),
\end{array}$$
for $j=1,2$.
Fix $j\ge 3$. Let $p_1\ge 2$ and $\frac 1{p_0}-\frac 1{p_1}=\frac 12$. By H\"older inequality,  for any $x\in B$, we have
$$\begin{array}{cl}
&\l(\dfrac
1{\Psi_\eta(B)|B|}\dint_{B(x_0,r)}|F(\sqrt{L})f_j|^{p_0}\,dy\r)^{1/p_0}\\
&\qquad\le C|B|^{-1/p_1}\Psi_\eta(B)^{-1/p_0}\|F^l(\sqrt{L})f_j\|_{L^{p_1}(\rz)}\\
&\qquad\le C|B|^{-1/p_1}\|F^l(\sqrt{L})\|_{L^{p_0}(U_j( B))\to L^{p_1}( B)}\Psi_\eta(B)^{-1/p_0}\|f_j\|_{L^{p_0}(\rz)}\\
&\qquad\le C 2^{j(\eta+n)/p_0}|B|^{\frac 12}\|F^l(\sqrt{L})\|_{L^{p_0}(U_j( B))\to L^{p_1}( B)}M_{\rho,\eta}(|f|^{p_0})^{1/p_0}(x).
\end{array}$$
Let $\frac 1{p_0}=\frac \tz 1+\frac {1-\tz}2$ and $\frac 1{p_1}=\frac \tz 2$, that is $\tz=2(\frac 1{p_0}-\frac 12)$. By interpolation,
$$\begin{array}{cl}
\|F^l(\sqrt{L})&\|_{L^{p_0}(U_j( B))\to L^{p_1}( B)}\\
&\le \|F^l(\sqrt{L})\|^{1-\tz}_{L^{2}(U_j( B))\to L^{\fz}( B)}\|\bar F^l(\sqrt{L})\|^{\tz}_{L^{2}(U_j( B))\to L^{\fz}( B)},
\end{array}$$
where $\bar F(L)$ be the operator with multiplier $\bar F$, the complex conjugate of $F$, and $\bar F$ satisfies the same estimates as $F$.

Next we estimate $\|F^l(\sqrt{L})\|_{L^{2}(U_j( B))\to L^{\fz}( B)}$. For every $l\in\zz$, let $K_{F^l}(\sqrt{L})(y,z)$ be the Schwartz kernel of operator $F^l(\sqrt{L})$. Then
$$\begin{array}{cl}
&\|F^l(\sqrt{L})\|^2_{L^{2}(U_j( B))\to L^{\fz}( B)}\\
&\qquad=\dsup_{y\in B}\dint_{U_j(B)}|K_{F^l}(\sqrt{L})(y,z)|^2dz\\
&\qquad \le C2^{-2sj}(2^lr)^{-2s}\dsup_{y\in B}\dint_\rz|K_{F^l}(\sqrt{L})(y,z)|^2(1+2^l|y-z|)^{2s}dz.
\end{array}\eqno(4.16)$$
Applying Lemma 4.3 with $F=F^l$ and $R=2^l$, we then have
$$\begin{array}{cl}
\dint_\rz|K_{F^l}(\sqrt{L})(y,z)|^2(1+2^l|y-z|)^{2s}dz&\le C 2^{ln}
\l(1 +\dfrac1{2^l\rho(y)}\r)^{-N}\|\dz_{2^l} (F^l)\|_{W^\fz_{s}}^2\\
&= C 2^{ln}
\l(1 +\dfrac r{\rho(y)}\dfrac 1{2^lr}\r)^{-N}\|\dz_{2^l} (F^l)\|_{W^\fz_{s}}^2\\
&\le C 2^{ln}\min\{1,(2^lr)^N\}\|\dz_{2^l} [\vz_l F]\|_{W^\fz_{s}}^2,
\end{array}\eqno(4.17)$$
since $r\ge \rho(x_0)$, so there a constant C such that $\rho(y)\le Cr$ for any $y\in B(x_0,r)$ by Lemma 2.1.

By (4.16) and (4.17), we get
$$\begin{array}{cl}
&\|F^l(\sqrt{L})\|_{L^{2}(U_j( B))\to L^{\fz}( B)}\\
&\le \l(2^{-2sj}2^{ln}(2^lr)^{-2s}\min\{1,(2^lr)^N\}\r)^{1/2}\|\dz_{2^l} [\vz_l F]\|_{W^\fz_{s}}.
\end{array}\eqno(4.18)$$
Now, we estimate $\|\bar F^l(\sqrt{L})\|_{L^{2}(U_j( B))\to L^{\fz}( B)}$. The calculations symmetric to  (4.16) with $\sup_{y\in B}$ replaced by $\sup_{z\in U_j(B)}$ yield. Then by Lemma 2.1, we have
$$\begin{array}{cl}
&\|\bar F^l(\sqrt{L})\|^2_{L^{2}(U_j( B))\to L^{\fz}( B)}\\
&\le C2^{-2sj}2^{ln}(2^lr)^{-2s}\dsup_{z\in U_j(B)}\l(1 +\dfrac r{\rho(z)}\dfrac 1{2^lr}\r)^{-N}\|\dz_{2^l} (F^l)\|_{W^\fz_{s}}^2\\
&\le C 2^{-2sj+jNl_0/(l_0+1)}2^{ln}(2^lr)^{-2s}\min\{1,(2^lr)^N\}\|\dz_{2^l} [\vz_l F]\|^2_{W^\fz_{s}}.
\end{array}\eqno(4.19)$$
Combining (4.18) and (4.19), we get
$$\begin{array}{cl}
&\l(\dfrac
1{\Psi_\eta(B)|B|}\dint_{B}|F(\sqrt{L})f_j|^{p_0}\,dy\r)^{1/p_0}\\
&\quad\le C2^{-js+j(\eta+n)/p_0+jNl_0/2(l_0+1)}(2^lr)^{-s}\l(\min\{1,(2^lr)^{N}\}\max\{1,(2^lr)^{\frac n2}\}\r)\\
&\qquad\qquad\times M_{\rho,\eta}(|f|^{p_0})^{1/p_0}(x)\dsup_{l\in\zz}\|\dz_{2^l}[\vz_lF]\|_{W^\fz_s}.
\end{array}\eqno(4.20)$$
Therefore,
$$\begin{array}{cl}
&\dsum_{j=3}^\fz\dsum_{l=-\fz}^\fz\l(\dfrac
1{\Psi_\eta(B)|B|}\dint_{B}|F^l(\sqrt{L})f_j|^{p_0}\,dy\r)^{1/p_0}\\
&\qquad\le C \dsum_{j=3}^\fz 2^{-js+j(\eta+n)/p_0+jNl_0/2(l_0+1)}\l(\dsum_{l=-\fz}^\fz(2^lr)^{-s}\min\{1,(2^lr)^{N}\}\max\{1,(2^lr)^{\frac n2}\}\r)\\
&\qquad\qquad\times M_{\rho,\eta}(|f|^{p_0})^{1/p_0}(x)\dsup_{l\in\zz}\|\dz_{2^l}[\vz_lF]\|_{W^\fz_s}\\
&\qquad\le C \dsum_{j=3}^\fz 2^{-js+j(\eta+n)/p_0+jNl_0/2(l_0+1)}\l(\dsum_{l: 2^lr>1}(2^lr)^{-s+\frac n2}+\dsum_{l: 2^lr\le 1}(2^lr)^{N-s}\r)\\
&\qquad\qquad\times M_{\rho,\eta}(|f|^{p_0})^{1/p_0}(x)\dsup_{l\in\zz}\|\dz_{2^l}[\vz_lF]\|_{W^\fz_s}\\
&\qquad\le CM_{\rho,\eta}(|f|^{p_0})^{1/p_0}(x)\dsup_{l\in\zz}\|\dz_{2^l}[\vz_lF]\|_{W^\fz_s},
\end{array}\eqno(4.21)$$
if $N>s>n+\eta+N/2$.

Combining estimates (4.15) and (4.21), we have proved (4.9), and then estimate (3.5) holds. Note that $T=F(L)$ and $A_{r_B}=I-(I-e^{-r_B^2L})^M$ commutes, it is easy to see that
$$\|A_{r_B}f\|_{L^\fz(c_0B)}\le CM_{\rho,\eta}(|f|^{p_0})^{\frac 1{p_0}}(x)\quad \forall\ x\in c_0B.$$
From this, (3.6) holds.  Thus, (4.3) is proved.

Let us turn to prove (4.4).\  We will adapt some similar arguments  in  \cite{m,dm}. We set $\wz\in A_1^{\rho,\eta}$ with  $\eta>0$. We know that $\wz\in A_1^{\rho,\eta}\subset A_p^{\rho,\eta}$ for every $1<p<\fz$. Then $F(\sqrt{L})$ is bounded on $L^p(\wz)$ for $1<p<\fz$. On the other hand, it is enough to prove the desired inequality for $f\in L_c^\fz$. If $\lz>0$, we
 consider the variant Calder\'on-Zygmund
decomposition of $f$ at level $\lz$ (ref \cite{cw}) and there exists a collection of balls $\{B_i\}$ such that $\{x\in\rz:\  M_{\rho,\eta}f(x)>\wt c\lz\}=\bigcup_i B_i$, where $\wt c$ is a positive constant depending only on $n$.

Now we decompose $f$ as $f=g+h=g+\sum_i h_i$, where
$$g(x)=f(x)\chi_{\rz\setminus \bigcup_i B_i}+\dsum_i\l(\dfrac 1{\Psi_\eta(B_i)|B_i|}\dint_{B_i}f(y)\rho_i(y)dy\r)\chi_{B_i}(x),$$
$$h_i(x)=f(x)\rho_i(x)-\l(\dfrac 1{\Psi_\eta(B_i)|B_i|}\dint_{B_i}f(y)\rho_i(y)dy\r)\chi_{B_i}(x),$$
$$\rho_i(x)=\dfrac{\chi_{B_i(x)}}{\sum_j\chi_{B_j(x)}}\chi_{\cup_j B_j}\quad {\rm and}\quad \chi_{\cup_j B_j}\le C.$$
From these and the definition of $ M_{\rho,\eta}f(x)$, we have
 $$\lz<\dfrac 1{\Psi_\eta(B_i)|B_i|}\dint_{B_i}|f(x)|dx\le C_\eta \lz,$$
If we denote $\oz=\bigcup_i B_i:=\bigcup_i B_i(x_i,r_i)$, then $|f(x)|\le\lz\
a.e.\ x\in\rz\setminus\oz$,
and $|g(x)|\le C_\eta\lz$ a.e. $x\in\rz$.
Observe that
$$\begin{array}{cl}
\wz\{x\in\rz: \ |F(\sqrt{L})f(x)|>\lz\}&\le\wz\{x\in\rz: \ |F(\sqrt{L})g(x)|>\lz/2\}\\
&\qquad+\wz\{x\in\rz: \ |F(\sqrt{L})h(x)|>\lz/2\}.
\end{array}$$
Since $F(\sqrt{L})$ is a bounded operator on $L^2(\wz)$. Then
$$\begin{array}{cl}
\wz\{x\in\rz: \ |F(\sqrt{L})g(x)|>\lz/2\}&\le \dfrac C{\lz^2}\dint_\rz |g(x)|^2\wz(x)dx\le \dfrac C{\lz}\dint_\rz |g(x)|\wz(x)dx\\
&\le \dfrac C{\lz}\l(\|f\|_{L^1(\wz)}+\dint_{\bigcup_iB_i}|g(x)|\wz(x)dx\r)\\
&\le \dfrac C{\lz}\|f\|_{L^1(\wz)}.
\end{array}$$
Let us deal with $F(\sqrt{L})h$. Set $t_i=r_i^2$ and write
$$h(x)=\dsum_i h_i(x)=\dsum_i(A_{t_i}h_i(x)+(h_i(x)-A_{t_i}h_i(x))),$$
where $A_{t_i}=e^{-t_i L}$.
Note that for any $N>2(l_0+1)\eta$, we have
$$\begin{array}{cl}
|A_{t_i}h_i(x)|&\le C_Nt_i^{-n/2}\dsup_{y\in B_i}e^{-b |x-y|^2/(t_i)}\dint_{B_i}\l(1+t_i/\rho^2(y)\r)^{-N}|b_i(y)|dy\\
&\le C_Nt_i^{-n/2}\l(1+t_i/\rho^2(x_i)\r)^{-N/(l_0+1)}\dint_{B_i}\l(1+t_i/\rho^2(y)\r)^{-N}|b_i(y)|dy\\
&\le C_Nt_i^{-n/2}\l(1+t_i/\rho^2(x_i)\r)^{-N/2(l_0+1)}\dsup_{y\in B_i}e^{-b|x-y|^2/(4t_i)}\Psi_\eta(B_i)^{-1}\dint_{B_i}|b_i(y)|dy\\
&\le C_N\lz|B_i|\l(1+t_i/\rho^2(x_i)\r)^{-N/2(l_0+1)}t_i^{-n/2}\dinf_{y\in B_i}e^{-b|x-y|^2/(4t_i)}\\
&\le C_N\lz\l(1+t_i/\rho^2(x_i)\r)^{-N/2(l_0+1)}t_i^{-n/2}\dint_\rz e^{-b|x-z|^2/(4t_i)}\chi_{B_i}(z)dz.
\end{array}$$
If $0\le u$ and  $\|u\|_{L^2(\rz)}=1$, by the weak type (1,1) of $M_{\rho,\eta}$, we get
$$
\dint_\rz \l(\dsum_i|A_{t_i}b_i(x)|\r)u(x)dx\le C\dint_\rz M_{\rho,\eta}u(y)\chi_{\bigcup_i B_i}(y)dy\le \lz^{1/2}\|f\|^{\frac 12}_{L^1(\rz)}.
\eqno(4.22)$$
From this, we know that
the above two series converge in $L^2(\rz)$.

In addition, it is easy to see that
$$\begin{array}{cl}
\wz\{x\in\rz: \ |F(\sqrt{L})h(x)|>\lz/2\}&\le \wz\{x\in\rz: \ |F(\sqrt{L})\l(\dsum_i(A_{t_i}h_i)\r)(x)|>\lz/4\}\\
&\ +\wz\{x\in\rz: \ |F(\sqrt{L})\l(\dsum_i(h_i-A_{t_i}h_i)\r)(x)|>\lz/4\}.
\end{array}$$
Since $T$ is bounded on $L^2(\wz)$, we have
$$\begin{array}{cl}
\wz\{x\in\rz: \ |F(\sqrt{L})\l(\dsum_i(A_{t_i}h_i)\r)(x)|>\lz/4\}&\le \dfrac 4{\lz^2}\dint_\rz |F(\sqrt{L})\l(\dsum_i(A_{t_i}h_i)\r)(x)|^2\wz(x)dz\\
&\le \dfrac C{\lz^2}\|\dsum_i|A_{t_i}h_i|^2\wz^{\frac 12}\|_{L^2(\rz)}^2.
\end{array}$$
We consider $0\le u\in L_c^\fz(\rz)$ with $\|u\|_{L^2}=1$. We apply (4.22) to $u\wz^{\frac 12}$, by the weighted weak type (1,1) of $M_{\rho,\eta}$(see Lemma 2.2), we then have
$$\begin{array}{cl}
\dint\l(\dsum_i|A_{t_i}h_i|^2\wz^{\frac 12}\r)u(x)dz&\le C\lz\dint_\rz M_{\rho,\eta}(u\wz^{\frac 12})(x)\chi_{\bigcup_i B_i}(x)dx\\
&\le C\lz \|M_{\rho,\eta}(u\wz^{\frac 12})\|_{L^2(\wz^{-1})}\|\chi_{\bigcup_i B_i}\|_{L^2(\wz)}\\
&\le C\lz^{\frac 12}\|f\|_{L^1(\wz)}.
\end{array}$$
Hence,
$$\wz\l\{x\in\rz: \ |F(\sqrt{L})\l(\dsum_iA_{t_i}h_i\r)(x)|>\lz/4\r\}\le \dfrac C{\lz}\|f\|_{L^1_\wz}.\eqno(4.23)$$
On the other hand, set $\wt B_i=2B_i$, we then have
$$\begin{array}{cl}
&\wz\l\{x\in\rz: \ |F(\sqrt{L})\l(\dsum_i(h_i-A_{t_i}h_i)\r)(x)|>\lz/4\r\}\\
&\qquad\le \wz\l( \bigcup_i\wt B_i\r)+\wz\l\{x\in\rz\setminus\bigcup_i\wt B_i: \ |F(L)\l(\dsum_i(h_i-A_{t_i}h_i)\r)(x)|>\lz/4\r\}\\
&\qquad\le \wz\l( \bigcup_i\wt B_i\r)+\dfrac 4\lz \dsum_i\dint_{\rz\setminus \wt B_i}|F(\sqrt{L})(h_i-A_{t_i}h_i)(x)|\wz(x)dx\\
&\qquad:=I_1+I_2.
\end{array}$$
For $I_1$, since $\wz\in A_1^{\rho,\eta}$, we get
$$\begin{array}{cl}
I_1&\le C\dsum_j\dfrac {\wz(\wt B_j)}{\Psi_\eta(\wt B_j)|\wt B_j|}\vz(B_j)|B_j|\\
&\le \dfrac C\lz\dsum_j\dfrac {\wz(\wt B_j)}{\Psi_\eta(\wt B_j)|\wt
B_j|}\dint_{B_j}|f(x)|dx\\
&\le \dfrac C\lz\dsum_j\dint_{B_j}|f(x)|M_{\rho,\eta}\wz(x)dx
\le \dfrac C\lz\dint_\rz|f(x)|\wz(x)dx.
\end{array}\eqno(4.24)$$
By Lemma 4.5,  we then have
$$\begin{array}{cl}
\dint_{\rz\setminus \wt B_i}&|F(\sqrt{L})(h_i-A_{t_i}h_i)(x)|\wz(x)dx\\
&\le
\dint_{B_i}|h_i(y)|\dint_{\rz\setminus B(y, r_i)}|K_{F(1-A_{t_i})(\sqrt{L})}(x,y)|\wz(x)dxdy\\
&\le C\dint_{B_i}|h_i(y)|\wz(y)dy\\
&\le C\dint_{B_i}|f(y)\wz(y)dy+C\lz\wz(B_i).
\end{array}\eqno(4.25)$$
Combining  (4.24) and (4.25), we get
$$I_2\le \dfrac C\lz\dint_\rz|f(x)|\wz(x)dx.$$
Thus, (4.4) holds.\hfill$\Box$

\subsection*{4.3. Proof of Theorem 4.2}
We first prove (4.5).  In fact,  adapting the similar proof of (4.3), we can prove (3.7), (3.8) and (3.9) hold, if $\dsup_{t>0}\|\eta\dz_t F\|_{W^\fz_{s}}\le C_s$ for $s$ large enough. We omit the details here.

It remains to prove (4.6). We borrow some ideas from \cite{pg,p1}. We consider the case $k=1$, the case $k\ge 2$ is similar. We give the same Calder\'on-Zygmund decomposition and $g$  and $h$ functions as in the proof of (4.4). Observe that $\dsup_{t>0}\|\eta \dz_t F\|_{W^\fz_s}\sim \dsup_{t>0}\|\eta \dz_t F\|_{W^\fz_s}$, where $G(\lz)=F(\sqrt{\lz})$. We write $F_b(\sqrt{L})f(x)=F_b^1(\sqrt{L})f(x)$.
Then,
$$\begin{array}{cl}
\wz\{x\in\rz: \ |F_b(\sqrt{L})f(x)|>\lz\}&\le\wz\{x\in\rz: \ |F_b(\sqrt{L})g(x)|>\lz/2\}\\
&\qquad+\wz\{x\in\rz: \ |F_b(\sqrt{L})h(x)|>\lz/2\}.
\end{array}$$
Since $F_b(\sqrt{L})$ is a bounded operator on $L^2(\wz)$. Then
$$\begin{array}{cl}
\wz\{x\in\rz: \ |F_b(\sqrt{L})g(x)|>\lz/2\}&\le \dfrac C{\lz^2}\dint_\rz |g(x)|^2\wz(x)dx\le \dfrac C{\lz}\dint_\rz |g(x)|\wz(x)dx\\
&\le \dfrac C{\lz}\l(\|f\|_{L^1_\wz}+\dint_{\bigcup_iB_i}|g(x)|\wz(x)dx\r)\le \dfrac C{\lz}\|f\|_{L^1_\wz}.
\end{array}$$
Let us deal with $F_b(\sqrt{L})h$. Set $t_i=r_i^2$ and write
$$h(x)=\dsum_i h_i(x)=\dsum_i(A_{t_i}h_i(x)+(h_i(x)-A_{t_i}h_i(x))).$$
Similar to the proof of (4.23), we have
$$\wz\l\{x\in\rz: \ |F_b(\sqrt{L})\l(\dsum_iA_{t_i}h_i\r)(x)|>\lz/4\r\}\le \dfrac C{\lz}\|f\|_{L^1(\wz)}.$$
In addition, set $\wt B_i=2B_i$ and $\wt \oz=\bigcup_i\wt B_i$, we then have
$$\begin{array}{cl}
&\wz\l\{x\in\rz: \ |F_b(\sqrt{L})\l(\dsum_i(h_i-A_{t_i}h_i)\r)(x)|>\lz/4\r\}\\
&\qquad\le \wz\l(\wt \oz\r)+\wz\l\{x\in\rz\setminus\wt \oz: \ |F_b(\sqrt{L})\l(\dsum_i(h_i-A_{t_i}h_i)\r)(x)|>\lz/4\r\}\\
&\qquad:=II_1+II_2.
\end{array}$$
For $II_1$, similar to the proof of (4.24), we get
$$
II_1\le \dfrac C\lz\dint_\rz|f(x)|\wz(x)dx.
$$
For $II_1$,  we have
$$\begin{array}{cl}
II_2
&\le \wz\l\{x\in\rz\setminus\wt \oz: \ |F_b(\sqrt{L})\l(\dsum_{i\in E_1}(h_i-A_{t_i}h_i)\r)(x)|>\lz/8\r\}\\
&\qquad+\wz\l\{x\in\rz\setminus\wt \oz: \ |F_b(\sqrt{L})\l(\dsum_{i\in E_2}(h_i-A_{t_i}h_i)\r)(x)|>\lz/8\r\}\\
&\le \wz\l\{x\in\rz\setminus\wt \oz: \ |\dsum_{i\in E_1}(b-b_{B_i})F(\sqrt{L})\l((h_i-A_{t_i}h_i)\r)(x)|>\lz/16\r\}\\
&\qquad+ \wz\l\{x\in\rz\setminus\wt \oz: \ |F(\sqrt{L})\l(\dsum_{i\in E_1}(b-b_{B_i})(h_i-A_{t_i}h_i)\r)(x)|>\lz/16\r\}\\
&\qquad+\wz\l\{x\in\rz\setminus\wt \oz: \ |F_b(\sqrt{L})\l(\dsum_{i\in E_2}(h_i-A_{t_i}h_i)\r)(x)|>\lz/8\r\}\\
&:=II_{21}+II_{22}+II_{23},
\end{array}$$
where $E_1=\{i:\ r_i<\rho(x_i)\}$ and $E_2=\{i:\ r_i\ge
\rho(x_i)\}$, $B_i=B(x_i,r_i)$ and $b_{B_i}$ is the average on $B_i$.

 To estimate $II_{21}$, by Lemma 4.7, we have
$$\begin{array}{cl}
II_{21}&\le \dfrac C{\lz}\dsum_{i\in E_1}\dint_{\rz\setminus
\wt\Omega}|b(x)-b_{B_i}||F(\sqrt{L})\l((h_i-A_{t_i}h_i)\r)(x)|\wz(x)dx\\
&\le \dfrac C{\lz}\dsum_{i\in E_1}\dint_{\rz\setminus
\wt\Omega}|b(x)-b_{B_i}|\dint_{B_i}|K_{F(I-A_{t_i})(\sqrt{L})}(x,y)||h_i(y)|dy\wz(x)dx\\
 &\le \dfrac C{\lz}\dsum_{i\in E_1}\dint_{B_i}|h_i(y)|\dint_{\rz\setminus
\wt B_i}|b(x)-b_{B_i}||K_{F(I-A_{t_i})(\sqrt{L})}(x,y)|\wz(x)dxdy\\
&\le \dfrac {C}{\lz}\dsum_{j\in E_1}\dint_{B_i}|h_i(y)|\wz(y)dy\le \dfrac {C}{\lz}\dint_\rz |f(x)|\wz(x)dx.
\end{array}$$
For $II_{22}$,
by (2.4) and (2.5), from Theorem 4.1 (b), we obtain
$$\begin{array}{cl}
II_{22}&\le \dfrac C{\lz}\dsum_{i\in E_1}\dint_{B_i}
|b(x)-b_{B_i}||h_j(x)|\wz(x)dx\\
&\le \dfrac C{\lz}\dsum_{i\in E_1}\dint_{B_i}
|b(x)-b_{B_i}||f(x)|\wz(x)dx\\
&\quad+\dfrac C{\lz}\dsum_{i\in E_1}\dfrac
1{|B_i|}\dint_{B_i}|f(y)|dy\dint_{B_i}
|b(x)-b_{B_i}|\wz(x)dx\\
&\le \dfrac C{\lz}\dsum_{i\in E_1}\wz(B_i)
\|b\|_{BMO_\tz(\rho)}\|f\|_{\log L, B_i,\wz}+\dfrac C{\lz}\dsum_{i\in E_1}\dfrac
{\wz(B_i)}{|B_i|}\dint_{B_i}|f(y)|dy
\|b\|_{BMO_\tz(\rho)}\\
&\le \dfrac C{\lz}\dsum_{i\in E_1}\wz(B_i)
\|b-b_{B_i}\|_{\exp L,B_i,\wz}\l(\lz+\dfrac
\lz{\wz(B_i)}\dint_{B_i}\Phi\l(\dfrac {|f(y)|}\lz\r)\wz(y)\,dy\r)\\
&\quad+\dfrac C{\lz}\dsum_{i\in E_1}\dfrac
{\wz(B_i)}{|B_i|}\dint_{B_i}|f(y)|dy
\|b\|_{BMO_\tz(\rho)}\\
&\le C\dint_\rz\Phi\l(\dfrac {|f(x)|}\lz\r)\wz(x)dx.
\end{array}$$
 For $II_{23}$. If $i\in E_2$, by Lemma 4.8, for any a  ball $B_i=B(x_i,r_i)$,  then the ball
$B$ can be decomposed into finite  balls $\{Q^j_i\}^{j_i}_{i=1}$, there exists a constant $L$ depending only on $n$ such that at most $L$ balls are disjoint each other, and $B\subset \bigcup_i^{m} Q_i\subset 2nB$ and
$r_i/2\le \rho(x)\le 2n C_0r_i$ for some $x\in
Q^j_i=Q(x^j_i,r^j_i)$. Let $h_i^j=h_i\chi_{Q_i^j\bigcap B_i}$, we have
$$\begin{array}{cl}
II_{23}
&\le \wz\{y\in\rz\setminus \wt\Omega:\ \sum_{j\in E_2}\dsum_{j=1}^{j_i}|b(x)-b_{Q^j_i}||F(\sqrt{L})((h^j_i-A_{t_i}h^j_i)(x)|>\lz/16\}\\
&\quad+\wz\{y\in\rz\setminus \wt\Omega:\ |F(\sqrt{L})\l(\sum_{i\in E_2}\dsum_{j=1}^{j_i}(b(x)-b_{Q^j_i})h^j_i\r)(x)|>\lz/16\}\\
&:=II_{23}^1+II_{23}^2.
\end{array}$$
Similar to the proof of $II^1_{23}$, note that  $r^j_i\sim
\rho(x^j_i)$, by Lemma 4.7, we then have
$$\begin{array}{cl}
II_{21}&\le \dfrac C{\lz}\dsum_{i\in
E_2}\dsum_{j=1}^{j_i}\dint_{\rz\setminus
\wt\Omega}|b(x)-b_{Q^j_i}||F(\sqrt{L})\l((h^j_i-A_{t_i}h^j_i)\r)(x)|\wz(x)dx\\
&\le \dfrac C{\lz}\dsum_{i\in
E_2}\dsum_{j=1}^{j_i}\dint_{\rz\setminus
\wt\Omega}|b(x)-b_{Q^j_i}|\dint_{Q^j_i}|K_{F(I-A_{t_i})(\sqrt{L})}(x,y)||h^j_i(y)|dy\wz(x)dx\\
 &\le \dfrac C{\lz}\dsum_{i\in
E_2}\dsum_{j=1}^{j_i}\dint_{Q^j_i}|h^j_i(y)|\dint_{\rz\setminus
2 Q^j_i}|b(x)-b_{Q^j_i}||K_{F(I-A_{t_i})(\sqrt{L})}(x,y)|\wz(x)dxdy\\
&\le \dfrac C{\lz}\dsum_{i\in
E_2}\dsum_{j=1}^{j_i}\dint_{Q^j_i\bigcap B_i}|h_i(y)|\wz(y)dy\le \dfrac {C}{\lz}\dsum_{i\in E_2}\dint_{B_i}|h_i(y)|\wz(y)dy\\
\end{array}$$
$$\begin{array}{cl}
&\le \dfrac {C}{\lz}\dsum_{i\in E_2}\dint_{B_i}|f(y)|\wz(y)dy+\dfrac {C}{\lz}\dsum_{i\in E_2}\dfrac
{\wz(B_i)}{\Psi_\eta(B_i)|B_i|}\dint_{B_i}|f(y)|dy\\
&\le \dfrac {C}{\lz}\dint_\rz |f(x)|\wz(x)dx.
\end{array}$$
Similar to the proof of $II_{22}$, set $f_i^j=f\chi_{Q_i^j\bigcap B_i}$, note that  $r^j_i\sim
\rho(x^j_i)$, we then have
$$\begin{array}{cl}
II^2_{23}&\le \dfrac C{\lz}\dsum_{i\in
E_2}\dsum_{j=1}^{j_i}\dint_{Q^j_i}
|b(x)-b_{Q^j_i}||h^j_i(x)|\wz(x)dx\\
&\le \dfrac C{\lz}\dsum_{i\in
E_2}\dsum_{j=1}^{j_i}\dint_{Q^j_i}
|b(x)-b_{Q^j_i}||f_i^j(x)|\wz(x)dx\\
&\quad+\dfrac C{\lz}\dsum_{i\in
E_2}\dsum_{j=1}^{j_i}\dfrac
1{\Psi_\eta(B_i)|B_i|}\dint_{B_i}|f(y)|dy\dint_{Q^j_i}
|b(x)-b_{Q^j_i}|\wz(x)dx\\
&\le \dfrac C{\lz}\dsum_{i\in
E_2}\dsum_{j=1}^{j_i}
\|b\|_{BMO_\tz(\rho)}\wz(Q^j_i)\|f_i^j\|_{\log L, Q^j_i,\wz}\\
&\quad+\dfrac C{\lz}\dsum_{i\in
E_2}\dsum_{j=1}^{j_i}\dfrac
{\wz(Q_i^j)}{\Psi_\eta(B_i)|B_i|}\dint_{B_i}|f(y)|dy
\|b\|_{BMO_\tz(\rho)}\\
&\le \dfrac C{\lz}\dsum_{i\in
E_2}\dsum_{j=1}^{j_i}\wz(Q^j_i)\l(\lz+\dfrac \lz{\wz(Q_i^j)}\dint_{Q_i^j\bigcap B_i}\Phi(|f(y)|/\lz)\wz(y)dy\r)\\
&\quad+\dfrac C{\lz}\dsum_{i\in
E_2}\dfrac
{L\wz(2n B_i)}{\Psi_\eta(B_i)|B_i|}\dint_{B_i}|f(y)|dy\\
&\le C\dsum_{i\in
E_2}\l(\wz(2n B_i)+\dint_{B_i}\Phi(|f(y)|/\lz)\wz(y)dy+\dfrac 1{\lz}\dint_{B_i}|f(y)|\wz(y)dy\r)\\
&\le C\dint_\rz\Phi\l(\dfrac {|f(x)|}\lz\r)\wz(x)dx.
\end{array}$$
In the last inequality we used the following fact that (see the proof of (4.24))
$$\dsum_{i\in
E_2}\wz(2n B_i)\le \dfrac {C}{\lz}\dint_\rz |f(x)|\wz(x)dx.$$
 Thus, (4.6) holds.\hfill$\Box$

\begin{center} {\bf 5. Littlewood-Paley operators}\end{center}
Let $L$ be the same as Section 4. We define the Littlewood-Paley operator for $x\in\rz$ and $f\in L^2(\rz)$,
$$g_Lf(x)=\l(\dint_0^\fz(tL)^{1/2}e^{-tL}f(x)|^2\dfrac {dt}t\r)^{1/2}.$$
We have the following $L^p$  estimates.

\begin{lem}\label{l5.1.}\hspace{-0.1cm}{\rm\bf 5.1.}\quad
Let $1<p<\fz$, then
$$\|g_Lf\|_{L^p}\sim \|f\|_{L^p}, \quad \forall f\in L^p\bigcap L^2.$$
\end{lem}
Lemma 5.1 is a special case in \cite{a}.

We have the following weighted  estimates for  Littlewood-Paley operators.
\begin{thm}\label{t5.1.}\hspace{-0.1cm}{\rm\bf 5.1.}\quad
Let $g_L$ be defined as above.
\begin{enumerate}
\item[(a)]  If  $1<p<\fz$ and $\wz\in A_p^{\rho,\fz}$, then
$$\|g_Lf\|_{L^p(\wz)}\le C \|f\|_{L^p(\wz)}, \quad \forall f\in L^\fz_c,$$
\item[(b)]  If  $1<p<\fz$ and $\wz\in A_p^{\rho,\fz}$, then
$$ \|f\|_{L^p(\wz)}\le C\|g_Lf\|_{L^p(\wz)}, \quad \forall f\in L^p\bigcap L^2,$$
\item[(c)]  If   $\wz\in A_1^{\rho,\fz}$, then
$$ \|g_Lf\|_{L^{1,\fz}(\wz)}\le C\|f\|_{{L^1(\wz)}}, \quad \forall f\in L^\fz_c.$$
\end{enumerate}
\end{thm}
We also define the commutator for  the Littlewood-Paley operator for $x\in\rz$ and $f\in L^2(\rz)$,
$$g_{L,b}^kf(x)=\l(\dint_0^\fz(tL)^{1/2}e^{-tL}(b(x)-b(\cdot))^kf(x)|^2\dfrac {dt}t\r)^{1/2}.$$
We give the following weighted  estimates the commutator for the Littlewood-Paley operator.
\begin{thm}\label{t5.2.}\hspace{-0.1cm}{\rm\bf 5.2.}\quad
Let $k\in\nn$ and $b\in BMO_\tz(\rho)$.
\begin{enumerate}
\item[(a)]  If  $1<p<\fz$ and $\wz\in A_p^{\rho,\fz}$, then
$$\|g^k_{L,b}f\|_{L^p(\wz)}\le C\|b\|^k_{BMO_\tz(\rho)} \|f\|_{L^p(\wz)},$$
\item[(b)]  If   $\wz\in A_1^{\rho,\fz}$, then
$$ \wz(\{x\in\rz: \ |g_{L,b}^kf(x)|>\lz\})\le C\Phi(\|b\|_{BMO_\tz(\rho)}) \dint_\rz\Phi\l(\dfrac
{|f(x)|}\lz\r)\wz(x)dx,$$ where
$\Phi(t)=t\log(e+t)^k$.
\end{enumerate}
\end{thm}
Before we begin the arguments, we recall some basic facts about Hilbert-valued extensions of scalar inequalities. To do so we introduce some notation: by ${\cal H}$ we mean $L^2((0,\fz),\frac {dt}t)$ and $\|\cdot\|$ denotes the norm in ${\cal H}$. Hence , for a function $h:\rz\times(0,\fz)\to\cc$, we have for $x\in\rz$
$$\|{\cal H}(x,\cdot)\|=\l(\dint_0^\fz|h(x,t)|^2\frac {dt}t\r)^{1/2}.$$
In particular,
$$g_Lf(x)=\|\vz(L,\cdot)f(x)\|$$
and
$$g_{L,b}^kf(x)=\|\vz(L,\cdot)(b(x)-b(\cdot))^kf(x)\|$$
with $\vz(z,t)=(tz)^{1/2}e^{-tz}$.
Let $L^p_{\cal H}(\wz)$ be the space of ${\cal H}$-valued $L^p(\wz)$-functions equipped with the
norm
$$\|h\|_{L_{\cal H}^p(\wz)}=\l(\dint_\rz\|h(x,\cdot)\|^pd\wz(x)\r)^{1/p}.$$
\begin{lem}\label{l5.2.}\hspace{-0.1cm}{\rm\bf 5.2(\cite{am}).}\quad
Let $\mu$ be a Borel measure on $\rz$. Let $1\le p\le q<\fz$. Let $D$ be a subspaces of ${\cal M}$, the space of measurable functions in $\rz$.
Let $S, T$ be linear operators from ${\cal D}$ into ${\cal M}$. Assume there exists $C_1>0$ such that for all $f\in{\cal D}$, we have
$$\|Tf\|_{L^q(\mu)}\le C_1\dsum_{j\ge 1}\az_j\|Sf\|_{L^q(F_j,\mu)},$$
where $F_j$ are subsets of $\rz$ and $\az_j\ge 0$. Then, there is an ${\cal H}$-valued extension with the same constant:
for all $f:\rz\times(0,\fz)\to\cc$ such that for (almost) all $t>0, f(\cdot,t)\in{\cal D}$,
$$\|Tf\|_{L^q_{\cal H}(\mu)}\le C_1\dsum_{j\ge 1}\az_j\|Sf\|_{L^q_{\cal H}(F_j,\mu)}.$$
\end{lem}
\subsection*{5.2. Proof of Theorem 5.1(a)}
Since $\wz\in A_p^{\rho,\fz}$, there exist $1<p_0<p$ and $\tz$ such that
$\wz\in A_{p/p_0}^{\rho,\tz}$. We are going to apply Theorem 3.2 with $k=0$, $T=g_L, A_r=I-(I-e^{-r^2L})^M$,
$M\in\nn$ large enough.  We first show (3.9) holds for all $f\in L_c^\fz$ and any $\eta>0$.

Let $1\le m\le M$. If $B=B(x_0,r)$  with $r<\rho(x_0)$, $\bar B=c_0B$, $B_j=2^{j+1}\bar B$, $C_j(\bar B)=B_j\setminus B_{j-1}$  for $j\ge 1$, and $g\in L^{p_0}$ with $\supp \, g\subset C_j(\bar B)$, by (4.1), we have
$$\|e^{-mr^2L}g\|_{L^\fz(\bar B)}\le C_12^{j(n+\eta)}e^{-\az 4^j}\l(\dfrac 1{\Psi_\eta(B_j)|B_j|}\dint_{C_j(\bar B)}|g|^{p_0}dx\r)^{1/p_0}.\eqno(5.1)$$
Lemma 5.2 applied to $S=I, \ T: L^{p_0}=L^{p_0}(\rz)\to L^{q_0}(\rz)$ given by
$$Tg=(C_12^{j(n+\tz)}e^{-\az 4^j})^{-1}\dfrac{|B_j|^{1/p_0}}{|\bar B|^{1/q_0}}\chi_{\bar B}e^{mr^2L}(\chi_{C_j(\bar B)}g)$$ yields for any $q_0$ satisfying $p_0<q_0<\fz$
$$\l(\dfrac 1{|B_j|}\dint_{\bar B}\|e^{-mr^2L}g(x,\cdot)\|^{q_0}dx\r)^{1/q_0}\le C_12^{j(n+\eta)}e^{-\az 4^j}\l(\dfrac 1{\Psi_\eta(B_j)|B_j|}\dint_{C_j(\bar B)}|g|^{p_0}dx\r)^{1/p_0}.\eqno(5.2)$$
Since (5.2) holds for any  $p_0<q_0<\fz$ , so
$$\|e^{-mr^2L}g(x,\cdot)\|_{L^\fz(\bar B)}\le C_12^{j(n+\eta)}e^{-\az 4^j}\l(\dfrac 1{\Psi_\eta(B_j)|B_j|}\dint_{C_j(\bar B)}|g|^{p_0}dx\r)^{1/p_0}.\eqno(5.3)$$
for all $g\in L_{\cal H}^{p_0}$ with $\supp\, g(x,\cdot)\subset C_j(\bar B)$ for each $t>0$ and some $\az>0$.

For $h\in L_{\cal H}^{p_0}$, we  write $h(x,t)=\dsum_{j\ge 1}h_j(x,t), \ x\in\rz,\ t>0,$
where $h_j(x,t)=$ $h(x,t)\chi_{C_j(B)}(x)$. Using (5.3), we have for $1\le m\le M$,
$$\begin{array}{cl}
\|\|e^{-mr^2L}h(x,\cdot)\|\|_{L^\fz(\bar B)}&\le \dsum_{j\ge 1} \|\|e^{-mr^2L}g(x,\cdot)\|\|_{L^\fz(\bar B)}\\
&\le C\dsum_{j\ge 1}2^{j(n+\eta)}e^{-\az 4^j}\l(\dfrac 1{\Psi_\eta(B_j)|B_j|}\dint_{C_j(\bar B)}|h|^{p_0}dx\r)^{1/p_0}.
\end{array}\eqno(5.4)$$
Take $h(x,t)=(tL)^{1/2}e^{-tL}f(x)$. Since $g_lf(x)=\|h(x,\cdot)\|$ and $f\in L_c^\fz$, $h\in L_{\cal H}^{p_0}$ by Lemma 5.1 and
$$g_L(e^{-mr^2L}f)(x)=\l(\dint_0^\fz|(tL)^{1/2}e^{-tL}e^{-mr^2L}f(x)|^2\dfrac {dt}t\r)^{1/2}=\|e^{-mr^2L}h(x,\cdot)\|.$$ Thus, (5.4) implies
$$\|g_L(e^{-mr^2L}f)\|_{L^\fz(\bar B)}
\le C\dsum_{j\ge 1}2^{j(n+\eta)}e^{-\az 4^j}\l(\dfrac 1{\Psi_\eta(B_j)|B_j|}\dint_{B_j}|g_L f|^{p_0}dx\r)^{1/p_0}.$$
and it follows that $g_L$ satisfies (3.9).

We now show (3.7) with $k=0$ holds for all $f\in L_c^\fz$ and any $\eta>0$. Set $B=B(x_0,r)$  with $r<\rho(x_0)$, $\bar B=c_0B$.
Write $f=\dsum_{j\ge 1}f_j$ as before. If $j=1$ we use that both $g_L$ and $(I-e^{-r^2L})^M$ are bounded on $L^{p_0}$, then
$$\l(\dfrac 1{|\bar B|}\dint_{\bar B}|g_L(I-e^{-r^2L})^Mf_1|^{p_0}dx\r)^{1/p_0}\le C
\l(\dfrac 1{|4\bar B|}\dint_{4\bar B}|f|^{p_0}dx\r)^{1/p_0}.\eqno(5.5)$$
For $j\ge 2$, we observe that
$$\begin{array}{cl}
g_L((I-e^{-r^2L})^Mf_j)(x)&=\l(\dint_0^\fz|(tL)^{1/2}e^{-tL}(I-e^{-r^2L})^Mf_j(x)|^2\dfrac {dt}t\r)^{1/2}\\
&=\|\vz(L,\cdot)f_j(x,\cdot)\|,\end{array}$$
where $\vz(z,t)=(tz)^{1/2}e^{-tz}(1-e^{-r^2z})^M$. Then $\vz(z,t)$ is a holomorphic function in $\sum_\mu=\{z\in \cc^*: |arg\ z|<\mu\}$ with $\mu\in (\nu,\pi)$, where $\nu\in [0,\pi/2)$ is defined by
$$\nu=\dsup\{|arg<Lf,f>|: \ f\in {\cal D}(L)\}.$$
Assume that $\nu<\tz<v<\mu<\pi/2$. As in \cite{a}, we then have
$$\vz(L,t)=\dint_{\Gamma_+}e^{-zL}\eta_+(z,t)dz+\dint_{\Gamma_-}e^{-zL}\eta_-(z,t)dz,$$
where $\Gamma_{\pm}$ is the half-ray $\rr^+e^{\pm i(\pi/2-\tz)}$,
$$\eta_{\pm}(z,t)=\dfrac 1{2\pi  i}\dint_{\pm}e^{\xi z}\vz(\xi,t)d\xi,\quad z\in \Gamma_{\pm},$$
where $\gz_{\pm}$ being the half-ray $\rr^+e^{\pm i\nu}$ and $\Gamma_{\pm}$ is the half-ray $\rr^+e^{\pm i(\pi/2-\tz)}$.
Note that
$$|\eta_{\pm}(z,t)|\le C\dfrac {t^{1/2}}{(|z|+t)^{3/2}}\dfrac{r^{2M}}{(|z|+t)^M},
\quad z\in\Gamma_{\pm}, t>0.$$
Thus,
$$\|\eta_{\pm}(z,\cdot)\|\le\l(\dint_0^\fz \dfrac {t^{1/2}}{(|z|+t)^{3/2}}\dfrac{r^{2M}}{(|z|+t)^M}\dfrac {dt}t\r)^{1/21}\le \dfrac {r^{2M}}{|z|^{M+1}}.$$
Applying Minkowski's inequality, by (4.1), we get
$$\begin{array}{cl}
&\l(\dfrac 1{|\bar B|}\dint_{\bar B}\l\|\dint_{\Gamma_+}e^{-zL}f_j(x)\eta_+(z,\cdot)dz\r\|^{p_0}dx\r)^{1/p_0}\\
&\qquad\le\l(\dfrac 1{|\bar B|}\dint_{\bar B}\l(\dint_{\Gamma_+}|e^{-zL}f_j(x)|\|\eta_+(z,\cdot)\||dz|\r)^{p_0}dx\r)^{1/p_0}\\
&\qquad\le\dint_{\Gamma_+}\l(\dfrac 1{|\bar B|}\dint_{\bar B}|e^{-zL}f_j(x)|^{p_0}dx\r)^{1/p_0}\dfrac{r^{2M}}{|z|^{M+1}}|dz|\\
&\qquad\le C 2^{j(n+\eta)}\dint_0^\fz \l(\frac r{\sqrt{s}}\r)^{n/2}e^{-\az\frac {4^jr^2}{s}}\dfrac {r^M}{s^{M+1}}ds\l(\dfrac 1{\Psi_\eta(B_j)|B_j|}\dint_{C_j(\bar B)}|f|^{p_0}dy\r)^{1/p_0}\\
&\qquad\le C2^{j(n+\eta-2M)}\l(\dfrac 1{\Psi_\eta (B_j)|B_j|}\dint_{B_j}|f|^{p_0}dy\r)^{1/p_0}
\end{array}$$
provided $2M>n+\eta$. This plus the corresponding term for $\Gamma_-$ yield
$$\l(\dfrac 1{|\bar B|}\dint_{\bar B}|g_L(I-e^{-r^2L})^Mf_j|^{p_0}dx\r)^{1/p_0}\le C 2^{j(n+\eta-2M)}\l(\dfrac 1{\Psi_\eta (B_j)|B_j|}\dint_{B_j}|f|^{p_0}dy\r)^{1/p_0}. \eqno(5.6)$$
Combining (5.5) and (5.6), we obtain (3.7) holds whenever $2M>n+\eta$.

We now show (3.8) with $k=0$ holds for all $f\in L_c^\fz$ and any $\eta>0$. Set $B=B(x_0,r)$  with $r\ge\rho(x_0)$.
Write $f=\dsum_{j\ge 1}f_j=\dsum_{j\ge 1}f\chi_{C_j(B)}$,  $B_j=2^{j+1}B$ and $C_j(B)=B_j\setminus B_{j-1}$.
$$\l(\dfrac 1{\Psi_\eta(B)| B|}\dint_{B}|g_L f_1|^{p_0}dx\r)^{1/p_0}\le C
\l(\dfrac 1{\Psi_\eta(4B)|4 B|}\dint_{4 B}|f|^{p_0}dx\r)^{1/p_0}.\eqno(5.7)$$
For $j\ge 2$, we observe that
$$
g_L(f_j)(x)=\l(\dint_0^\fz|(tL)^{1/2}e^{-tL}f_j(x)|^2\dfrac {dt}t\r)^{1/2}
=\|\vz(L,\cdot)f_j(x,\cdot)\|,$$
where $\vz(z,t)=(tz)^{1/2}e^{tz}$. The functions $\eta_{\pm}(\cdot,t)$ associated  with $\vz(\cdot,t)$ verify
$$|\eta_{\pm}(z,t)|\le C\dfrac {t^{1/2}}{(|z|+t)^{3/2}},
\quad z\in\Gamma_{\pm}, t>0.$$
From this,  note that $r\ge\rho(x_0)$, By (4.1) and Lemma 2.1, we then have
$$\begin{array}{cl}
&\l(\dfrac 1{\Psi_\eta(B)|B|}\dint_{ B}\l\|\dint_{\Gamma_+}e^{-zL}f_j(x)\eta_+(z,\cdot)dz\r\|^{p_0}dx\r)^{1/p_0}\\
&\quad\le\l(\dfrac 1{\Psi_\eta(B)| B|}\dint_{ B}\l(\dint_{\Gamma_+}|e^{-zL}f_j(x)|\|\eta_+(z,\cdot)\||dz|\r)^{p_0}dx\r)^{1/p_0}\\
&\quad\le\dint_{\Gamma_+}\l(\dfrac 1{\Psi_\eta(B)| B|}\dint_{ B}|e^{-zL}f_j(x)|^{p_0}dx\r)^{1/p_0}\dfrac 1{|z|}|dz|\\
&\quad\le C 2^{j(n+\eta)}\dint_0^\fz \l(\frac 1{\sqrt{s}}\r)^{n}e^{-\az\frac {4^jr^2}{s}}\dint_B\l(1+\frac{\sqrt{s}}{\rho(x)}\r)^{-N}dx\dfrac 1sds\l(\dfrac 1{\Psi_\eta(B_j)|B_j|}\dint_{C_j(B)}|f|^{p_0}dy\r)^{1/p_0}\\
\end{array}$$
$$\begin{array}{cl}
&\quad\le C 2^{j(n+\eta)}\dint_0^\fz \l(\frac 1{\sqrt{s}}\r)^{n}e^{-\az\frac {4^jr^2}{s}}\dint_B\l(1+\l(\frac r{\rho(x_0)}\r)^{l_0/(l_0+1)} \frac {\sqrt{s}}r\r)^{-N}dx\dfrac 1sds\\
&\qquad\qquad\times\l(\dfrac 1{\Psi_\eta(B_j)|B_j|}\dint_{C_j(B)}|f|^{p_0}dy\r)^{1/p_0}\\
&\quad\le C 2^{j(n+\eta)}\dint_0^\fz \l(\frac r{\sqrt{s}}\r)^{n}e^{-\az\frac {4^jr^2}{s}}\min\{(\frac r{\sqrt{s}})^{N},1\}\dfrac 1sds\l(\dfrac 1{\Psi_\eta(B_j)|B_j|}\dint_{C_j(B)}|f|^{p_0}dy\r)^{1/p_0}\\
&\quad\le C 2^{j(n+\eta)}\l(\dint_0^{r^2} \l(\dfrac r{\sqrt{s}}\r)^{n}\l(\dfrac s{4^jr^2}\r)^{N}\dfrac 1sds+\dint_{r^2}^\fz \l(\frac r{\sqrt{s}}\r)^{n+N}\l(\dfrac s{4^jr^2}\r)^{N/2}\dfrac 1sds\r)\\
&\qquad\qquad\times\l(\dfrac 1{\Psi_\eta(B_j)|B_j|}\dint_{C_j(B)}|f|^{p_0}dy\r)^{1/p_0}\\
&\quad\le C2^{j(n+\eta-N)}\l(\dfrac 1{\Psi_\eta(B_j)|B_j|}\dint_{B_j}|f|^{p_0}dy\r)^{1/p_0}
\end{array}$$
provided $N>n+\eta$. This plus the corresponding term for $\Gamma_-$ yield
$$\l(\dfrac 1{\Psi_\eta(B)| B|}\dint_{ B}|g_Lf|^{p_0}dx\r)^{1/p_0}\le C 2^{j(n+\eta-N)}\l(\dfrac 1{\Psi_\eta(B_j)|B_j|}\dint_{B_j}|f|^{p_0}dy\r)^{1/p_0}.\eqno(5.8) $$
Combining (5.7) and (5.8), we obtain (3.8) holds whenever $N>n+\eta$.\hfill$\Box$

\subsection*{5.3. Proof of Theorem 5.1(b)}
To prove Theorem 5.1 (b), we introduce the following operator. Define for $f\in L_{\cal H}^2$ and $x\in\rz$,
$$T_Lf(x)=\dint_0^\fz(tL)^{1/2}e^{-tL}f(x,t)\dfrac {dt}t.$$
Recall that $(tL)^{1/2}e^{-tL}f(x,t)=(tL)^{1/2}e^{-tL}(f(\cdot,t))(x)$. Hence, $T_L$ maps ${\cal H}$-valued functions to $\cc$-valued
functions. For $f\in L_{\cal H}^2$ and $h\in L^2$, we have
$$\dint_\rz T_Lf\bar hdx=\dint_\rz\dint_0^\fz f(x,t)\overline{(tL^*)^{1/2}e^{-tL^*}h(x)}\dfrac {dt}tdx,$$
where $L^*$ is the adjoint of $L$, hence,
$$\l|\dint_\rz T_L f\bar hdx\r|\le \dint_\rz\|f(x,\cdot)\|g_{L^*}(h)(x)dx.$$
Since $g_{L^*}$ is bounded on $L^p$ for $1<p<\fz$. This and a density imply that $T_L$ has a bounded extension from $L_{\cal H}^p$ to $L^p$.
We next give a weighted result for the operator $T_L$.
\begin{lem}\label{l5.3.}\hspace{-0.1cm}{\rm\bf 5.3.}\quad
Let $1<p<\fz$ and $\wz\in A_p^{\rho,\fz}$, then for all $f\in L_c^\fz(\rz\times(0,\fz))$ we have
$$\|T_Lf\|_{L^p(\wz)}\le C\|f\|_{L^p(\wz)}.$$
Hence, $T_L$  has a bounded extension from $L^p_{\cal H}(\wz)$ to $L^p(\wz)$.
\end{lem}
{\it Proof.}\quad  We will apply Theorem 3.2 with $k=0$(its vector-valued extension) with underlying measure $dx$ and weight $\wz$ to linear operator $T=T_L$ and $A_r=I-(I-e^{-r^2L})^M$, $M\in\nn$ large enough. Adapting similar to the proof of Theorem 5.1, we can obtain the desired result. We omit the details here.\hfill$\Box$

\medskip

{\it Proof of Theorem 5.1 (b).}\quad Let $f\in L^2$ and define $F$  by $F(x,t)=(tL)^{1/2}e^{-tL}f(x)$. Note that $F\in L_{\cal H}^2$ since $\|F\|_{L_{\cal H}^2}=\|g_Lf\|_{L^2}$. By functional calculus on $L^2$, we have
$$f=2\dint_0^\fz(tL)^{1/2}e^{-tL}F(\cdot,t)\dfrac {dt}t=2T_L F\eqno(5.11)$$
with convergence in $L^2$. Note that $e^{-tL}$ has an infinitesimal generator on $L^p(\wz)$ for $1<p<\fz$. Let us call $L_{p,\wz}$ this generator. In particular $e^{-tL}$ and $e^{-tL_{p,\wz}}$ agree on $L^p(\wz)\bigcap L^2$. Our assert that $L_{p,\wz}$ has a bounded holomorphic functional calculus on $L^p(\wz)$, hence replacing $L$ by $L_{p,\wz}$ and $f\in L^2$ by $f\in L^p(\wz)$, we see that $F\in L_{\cal H}^p(\wz)$ with
$\|F\|_{L^p(\wz)}=\|g_{L_{P,\wz}}f\|_{L^p(\wz)}$ and (5.11) is valid with convergence in $L^p(\wz)$. Thus, by Lemma 5.3,
$$\|f\|_{L^p(\wz)}=2\|T_{L_{p,\wz}}F\|_{L^p(\wz)}\le C\|F\|_{L_{\cal H}^p(\wz)}=\|g_{L_{p,\wz}}f\|_{L^p(\wz)}.$$ Note that $g_L f=g_{L_{p,\wz}}f$ when $f\in L^2\bigcap L^p(\wz)$ and $T_L=T_{L_{p,\wz}}F$ when $F\in L_{\cal H}^2\bigcap L_{\cal H}^p(\wz)$.\hfill$\Box$
\subsection*{5.4. Proof of Theorem 5.1(c)}
To prove Theorem 5.1(c), it suffices to show the following Lemma.
\begin{lem}\label{l5.4.}\hspace{-0.1cm}{\rm\bf 5.4.}\quad
Suppose $\wz\in A^{\rho,\tz}_1$ and $B=B(x_0,r)$. Then for $\supp f\subset B$,
$$\dint_{\rz\setminus 2B(x_0,r)}|g_L(1-e^{-r^2L})f(x)|\wz(x)dx\le C\dint_B |f(y)|\wz(y)dy.\eqno(5.9)$$
\end{lem}
{\it Proof.}\quad In fact,
$$\begin{array}{cl}
\dint_{\rz\setminus 2B(x_0,r)}|g_L(1-e^{-r^2L})f(x)|\wz(x)dx&\le
\dint_{\rz\setminus 2B(x_0,r)}\l\|\dint_{\Gamma_+}e^{-zL}f(x)\eta_+(z,\cdot)dz\r\|\wz(x)dx\\
&\ +\dint_{\rz\setminus 2B(x_0,r)}\l\|\dint_{\Gamma_-}e^{-zL}f(x)\eta_+(z,\cdot)dz\r\|\wz(x)dx\\
&:=I_1+I_2.
\end{array}$$
We only give the estimate for $I_1$, $I_2$ is similar. Then for any $N$ large enough, we have
$$\begin{array}{cl}
I_1&\le \dint_{\rz\setminus 2B(x_0,r)}\dint_{\Gamma_+}|e^{-zL}f(x)|\|\eta_+(z,\cdot)\||dz|\wz(x)dx\\
&\le \dint_{\Gamma_+}\dint_{\rz\setminus 2B(x_0,r)}|e^{-zL}f(x)|\wz(x)dx\dfrac{r^{2}}{|z|^{2}}|dz|\\
&\le C\dint_0^\fz\dint_{\rz\setminus 2B(x_0,r)}\dint_B\l(1+\dfrac {\sqrt{s}}{\rho(y)}\r)^{-N}s^{-n/2}e^{-\az|x-y|^2/2s}|f(y)|dy\wz(x)dx\dfrac{r^{2}}{s^{2}}ds\\
&\le C\dint_B\dint_0^\fz \dint_{\rz\setminus B(y,r)}s^{-n/2}\l(1+\dfrac {\sqrt{s}}{\rho(y)}\r)^{-N}e^{-\az|x-y|^2/2s}\wz(x)dx e^{-\az r^2/2s}\dfrac{r^{2}}{s^{2}}ds|f(y)|dy\\
&\le C\dint_B\dint_0^\fz e^{-\az r^2/2s}\dfrac{r^{2}}{s^{2}}ds|f(y)|\wz(y)dy\le C\dint_B|f(y)|\wz(y)dy.
\end{array}$$
Thus, (5.9) holds. \hfill$\Box$
\subsection*{5.5. Proof of Theorem 5.2}
We first prove Theorem 5.2(a).  In fact,  adapting the similar proof of Theorem 5.1(a), we can prove (3.7), (3.8) and (3.9) hold. We omit the details here.

To prove Theorem 5.2(b), it suffices to show the following Lemma.
\begin{lem}\label{l5.5.}\hspace{-0.1cm}{\rm\bf 5.5.}\quad
Suppose $b\in BMO_\tz(\rho)$, $\wz\in A^{\rho,\tz_1}_1$,  and $B=B(x_0,r)$ with $r<\rho(x_0)$. Then for any $f\in L^1(\wz)$ and $\supp f\subset B$
$$\dint_{\rz\setminus 2B(x_0,r)}|(b(x)-b_B)g_L(1-e^{-r^2L})f(x)|\wz(x)dx\le C\dint_B |f(y)|\wz(y)dy.\eqno(5.10)$$
\end{lem}
{\it Proof.}\quad In fact,
$$\begin{array}{cl}
\dint_{\rz\setminus 2B(x_0,r)}&|(b(x)-b_B)g_L(1-e^{-r^2L})f(x)|\wz(x)dx\\&\le
\dint_{\rz\setminus 2B(x_0,r)}\l\|\dint_{\Gamma_+}(b(x)-b_B)e^{-zL}f(x)\eta_+(z,\cdot)dz\r\|\wz(x)dx\\
&\quad+\dint_{\rz\setminus 2B(x_0,r)}\l\|\dint_{\Gamma_-}(b(x)-b_B)e^{-zL}f(x)\eta_+(z,\cdot)dz\r\|\wz(x)dx\\
&:=I_1+I_2.
\end{array}$$
We only give the estimate for $I_1$, $I_2$ is similar. Then for any $N$ large enough, we have
$$\begin{array}{cl}
I_1&\le \dint_{\rz\setminus 2B(x_0,r)}\dint_{\Gamma_+}|e^{-zL}f(x)|\|\eta_+(z,\cdot)\||dz|\wz(x)dx\\
&\le \dint_{\Gamma_+}\dint_{\rz\setminus 2B(x_0,r)}|e^{-zL}f(x)|\wz(x)dx\dfrac{r^{2}}{|z|^{2}}|dz|\\
&\le C\dint_0^\fz\dint_{\rz\setminus 2B(x_0,r)}\dint_B\l(1+\dfrac {\sqrt{s}}{\rho(y)}\r)^{-N}s^{-n/2}(b(x)-b_B)e^{-\az|x-y|^2/2s}\\
&\qquad\times|f(y)|dy\wz(x)dx\dfrac{r^{2}}{s^2}ds\\
&\le C\dint_B\dint_0^\fz \dint_{\rz\setminus B(y,r)}s^{-n/2}\l(1+\dfrac {\sqrt{s}}{\rho(y)}\r)^{-N}(b(x)-b_B)e^{-\az|x-y|^2/2s}\wz(x)dx\\
&\qquad\times e^{-\az r^2/2s}\dfrac{r^{2}}{s^{2}}ds|f(y)|dy\\
&\le C\|b\|_{BMO_\tz(\rho)}\dint_B\dint_0^\fz e^{-\az r^2/2s}\dfrac{r^{2}}{s^{2}}ds|f(y)|\wz(y)dy\\
&\le C\|b\|_{BMO_\tz(\rho)}\dint_B|f(y)|\wz(y)dy.
\end{array}$$
Thus, (5.10) holds. \hfill$\Box$

 LMAM, School of Mathematical  Science,
 Peking University,
 Beijing, 100871,
People¡¯s Republic of China

 E-mail address:  tanglin@math.pku.edu.cn


\begin{center} {\bf References}\end{center}
\begin{enumerate}
\vspace{-0.3cm}
\bibitem[1]{a}P. Auscher,
On necessary and sufficient conditions for $L^p$-estimates of Riesz transforms associated to elliptic operators on $\rz$
and related estimates. Mem. Amer. Math. Soc. 186 (2007), no. 871, xviii+75 pp.
\vspace{-0.3cm}
\bibitem[2]{aa}P. Auscher and B. Ali,
Maximal inequalities and Riesz transform estimates on $L^p$ spaces for Schr\"odinger operators with nonnegative potentials. Ann. Inst. Fourier (Grenoble) 57 (2007), 1975-2013.
\vspace{-0.3cm}
\bibitem[3]{am} P. Auscher and J. Martell,
Weighted norm inequalities, off-diagonal estimates and elliptic operators.
I. General operator theory and weights, Adv. Math. 212 (2007), 225-276.
\vspace{-0.3cm}
\bibitem[4]{am1} P. Auscher and J. Martell,
Weighted norm inequalities, off-diagonal estimates and elliptic operators. III. Harmonic analysis
of elliptic operators. J. Funct. Anal. 241 (2006),  703-746.
\vspace{-0.3cm}
\bibitem[5]{bhs1} B. Bongioanni, E. Harboure and O. Salinas,
Commutators of Riesz transforms related to Schr\"odinger operators, J. Fourier Ana Appl. 17(2011), 115-134.
\vspace{-0.3cm}
\bibitem[6]{bhs2} B. Bongioanni, E. Harboure and O. Salinas,
Class of weights related to Schr\"odinger operators, J. Math. Anal. Appl. 373(2011), 563-579.
\vspace{-0.3cm}
\bibitem[7]{cd} F. Cacciafesta and P. D'Ancona,
Weighted $L^p$ estimates for powers of selfajoint operators, Advance in Math. 229(2012), 501-530.
\vspace{-0.3cm}
\bibitem[8]{cw}
R. Coifman and R. G.Weiss, Analyse harmonique non-commutative sur
certains espaces homog\`{e}nes. Lecture Notes in Mathematics, 242.
Springer, Berlin-New York, 1971.
\vspace{-0.3cm}
\bibitem[9]{d}
E. Davies, Uniformly elliptic operators with measurable coefficients. J. Funct. Anal. 132 (1995),  141-169.
\vspace{-0.3cm}
\bibitem[10]{dm}
X. T. Duong and A. MacIntosh,  Singular integral operators with non-smooth kernels on irregular domains.
Rev. Mat. Iberoamericana 15 (1999), 233-265.
\vspace{-0.3cm}
\bibitem[11]{dos}
X. T. Duong, M. Ouhabaz and A. Sikora,  Plancherel-type estimates and sharp spectral multipliers.
J. Funct. Anal. 196 (2002), 443-485.
\vspace{-0.3cm}
\bibitem[12]{dsy} X. T. Duong, A. Sikora and  L. Yan,
 Weighted norm inequalities, Gaussian bounds
and sharp spectral multipliers, J. Funct. Anal. 260 (2011),  1106-1131.
  \vspace{-0.3cm}
\bibitem[13]{dz}J. Dziuba\'{n}ski,
A spectral multiplier theorem for $H^1$ spaces associated with Schr\"o-dinger operators with potentials satisfying a reverse H\"older inequality, Illinois J. Math. 45(2001), 1301-1313.
\vspace{-0.3cm}
\bibitem[14]{dz1}J. Dziuba\'{n}ski,
Note on $H^1$ spaces related to degenerate Schr\"odinger operators, Illinois J. Math
49(2005), 1271-1297.
\vspace{-0.3cm}
\bibitem[15]{gr} J. Garc\'ia-Cuerva and J. Rubio de Francia,
Weighted norm inequalities and related topics, Amsterdam- New
York, North-Holland, 1985.
\vspace{-0.3cm}
\bibitem[16]{h} W. Hebish,
A multiplier theorem for Schr\"odinger operators, Collq. Math. 60/61 (1990), 659-664.
\vspace{-0.3cm}
\bibitem[17]{m}J. Martell,
Sharp maximal functions associated with approximations of the identity in spaces of homogeneous type and applications. Studia Math. 161(2004),113-145.
 \vspace{-0.3cm}
\bibitem[18]{p1}C. P\'erez,
Endpoint estimates for commutators of singular integral operators,
J. Funct. Anal. 128(1995), 163-185.
\vspace{-0.3cm}
\bibitem[19]{pg}C. P\'erez and R. Gonz\'alez,
Sharp weighted estimates for vector-valued singular integral
operators and commutators, Tohoku Math. J. 55(2003), 109-129.
\vspace{-0.3cm}
\bibitem[20]{rr} M. M. Rao and Z. D. Ren,
Theory of Orlicz spaces, Monogr. Textbooks Pure Appl. Math.146,
Marcel Dekker, Inc., New York, 1991.
\vspace{-0.3cm}
\bibitem[21]{s1}Z. Shen,
$L^p$ estimates for Schr\"odinger operators with certain
potentials, Ann. Inst. Fourier. Grenoble, 45(1995), 513-546.
\vspace{-0.3cm}
\bibitem[22]{s}  E. M. Stein,
Harmonic Analysis: Real-variable Methods, Orthogonality, and
Oscillatory integrals. Princeton Univ Press. Princeton, N. J.
1993.
\vspace{-0.3cm}
\bibitem[23]{t1} L. Tang,
Weighted norm inequalities for  Schr\"odinger type operators, arXiv: 1109. 0099.
\vspace{-0.3cm}
\bibitem[24]{t2} L. Tang,
Weighted norm inequalities for commutators of Littlewood-Paley functions related to Schr\"odinger operators, arXiv: 1109.0100.
\vspace{-0.3cm}
\bibitem[25]{yz}D. Yang, Y. Zhou,
Localized Hardy spaces $H^1$ related to admissible functions on RD-spaces and applications to Schr\"odinger operators,
Trans. Amer. Math. Soc.  363 (2011), 1197-1239.
\vspace{-0.3cm}
\bibitem[26]{z}  J. Zhong,
Harmonic analysis for some Schr\"odinger type operators, Ph.D. Thesis. Princeton University,
  1993.

\end{enumerate}
\end{document}